\documentclass[12pt]{article}
\usepackage{amsmath, exscale, epsfig,amssymb, 
 makeidx}

\small
\normalsize
\usepackage[latin1]{inputenc}
\usepackage[T1]{fontenc}
\usepackage{amsfonts}
\usepackage{amsbsy}
\usepackage{amscd}
\usepackage{theorem}

\usepackage{epsf}
\usepackage{psfrag}
\usepackage{graphicx}

\usepackage{multicol}
\allowdisplaybreaks

\def\build#1_#2^#3{\mathrel{\mathop{\kern 0pt#1}\limits_{#2}^{#3}}}
\def\noi{{\noindent}}

\def\cq{$\hfill \square$}
\def\un{{\bf 1}}

\newcommand{\bbE}{\mathbb{E}}
\newcommand{\bE}{{\bf E}}

\newcommand{\bN}{\mathbb{N}}
\newcommand{\bbN}{\mathbb{N}}
\newcommand{\bbP}{\mathbb{P}}

\newcommand{\bP}{{\bf P}}
\newcommand{\bM}{{\bf M}}
\newcommand{\bR}{\mathbb{R}}
\newcommand{\bZ}{\mathbb{Z}}

\newcommand{\cW}{\mathcal{W}}
\newcommand{\cD}{\mathcal{D}}

\newcommand{\cF}{\mathcal{F}}
\newcommand{\cQ}{\mathcal{Q}}
\newcommand{\cG}{\mathcal{G}}
\newcommand{\cJ}{\mathcal{J}}

\newcommand{\cM}{\mathcal{M}}
\newcommand{\cN}{\mathcal{N}}
\newcommand{\cP}{\mathcal{P}}
\newcommand{\cR}{\mathcal{R}}

\newcommand{\cT}{\mathcal{T}}

\def\w{{\rm w}}

\def\be{\begin{equation}}
\def\ee{\end{equation}}
\def\ba{\begin{eqnarray*}}
\def\ea{\end{eqnarray*}}

\def\noi{\noindent}

\def\supp{{\rm supp\,}}

\def\cqfd{ \hfill $\blacksquare$ }

\newtheorem{theorem}{Theorem}[section]
\newtheorem{lemma}[theorem]{Lemma}

{\theorembodyfont{\rmfamily}}
{\theorembodyfont{\rmfamily}\newtheorem{remark}{Remark}[section]}
\begin{document}

\title{ {\bf THE  PACKING MEASURE OF THE RANGE OF SUPER-BROWNIAN MOTION}.}
\author{Thomas {\sc Duquesne}  
\thanks{Laboratoire de Probabilit\'es et Mod\`eles Al\'eatoires; 
Universit\'e Paris 6, 16 rue Clisson, 75013 PARIS, FRANCE. Email: thomas.duquesne@upmc.fr } }

\vspace{2mm}
\date{\today} 

\maketitle

\begin{abstract} We prove that the total range of Super-Brownian motion with quadratic branching mechanism has an exact packing measure 
with respect to the gauge function $g(r)=r^4 (\log \log1/r)^{-3}$ in super-critical dimensions $d\geq 5$.  More precisely, we prove that the total occupation measure of Super-Brownian motion is equal to the $g$-packing measure restricted to its range, up to a deterministic multiplicative constant that only depends on space dimension $d$. \\

\end{abstract}

\section{Introduction}
\label{introd}
$\; $ 
The purpose of this paper is to provide an exact packing gauge function for the range of Super-Brownian motion with quadratic branching mechanism in super-critical dimensions $d \geq 5$. Dawson, Iscoe and Perkins \cite{DaIsPe} have proved that $h(r)=r^4 \log \log (1/r)$ is the exact Hausdorff gauge function for the range of Super-Brownian motion in dimensions $d\geq 5$; Le Gall \cite{LGRange} showed that $h(r)=   r^4 \log (1/r) \log \log \log (1/r)$ is the correct Hausdorff gauge function in critical dimension $d= 4$; by use of  Brownian Snake techniques, he proves that the total occupation measure of Super-Brownian motion is equal to the $h$-Hausdorff measure restricted to its range (up to an unknown deterministic multiplicative constant). Similarly, we prove in this paper that in dimensions $d\geq 5$, the total occupation measure of Super-Brownian motion coincides with the $g$-packing measure in $\bR^d$ restricted to its range, where $g(r)= r^4 (\log \log 1/r)^{-3}$.  This result contrasts with known results concerning the support of Super-Brownian motion at a fixed time: Le Gall, Perkins and Taylor \cite{LGPerTay95} prove that in dimensions $d \geq 3$ there is no 
exact packing function for the support of Super-Brownian motion and they also provide an optimal test in dimensions $d \geq 3$ (and a partial result in the critical dimension $d= 2$).

   Let us mention that the results of our paper apply to Integrated Super-Brownian Excursion measure (ISE) that is the scaling limit of various models in statistical mechanics in high dimensions (see Slade \cite{Slade} for a survey on this topic). 

\medskip

  Let us briefly state our main results: denote by $M_f(\bR^d)$ the set of finite measures defined on the Borel sets of $\bR^d$ equipped with the topology of weak convergence. For any $\mu \in  M_f (\bR^d)$, we denote by $\supp \mu$ its topological support that is the smallest closed subset supporting $\mu$. Super-Brownian motion with quadratic branching mechanism 
is a time-homogeneous $M_f(\bR^d)$-valued Markov process $(Z_t,t\geq 0; \bP_\mu , \mu \in M_f(\bR^d))$ whose transition kernels are characterized as follows: for every $\mu  \in M_f(\bR^d)$ and
for any continuous nonnegative function $f$, we have
$$\bE_\mu   \left[  \exp \left(-\langle Z_t, f \rangle \right) \right]
=\exp(-\langle \mu , u_t \rangle),$$
where the function $(u_t(x);t\geq 0,x\in \bR^d)$ is the unique nonnegative
solution of the integral equation
$$u_t(x)+  2\beta \int_0^t   K_{t-s} (u^2_s) (x)  \, ds = K_t (f) (x) \; , \quad x \in \bR^d \; , \; t \in [0, \infty).   $$
Here $\beta $ is a positive constant determining the branching rate and $(K_t , t \geq 0)$ stands for the transition semi-group of the standard $d$-dimensional Brownian motion.  We refer to Dynkin \cite{CRM}, Le Gall \cite{LG99} and Perkins \cite{Pe} for a general introduction on super-processes.

Let us fix $\mu \in M_f (\bR^d)$ and let us consider $Z=(Z_t, t \geq 0)$ under $\bP_\mu$. We assume that $Z$ is cadlag. We define the total range of $Z$ by 
\begin{equation}
\label{defboldR}
 {\bf R} = \bigcup_{\varepsilon >0} \overline{ \bigcup_{ a \geq \varepsilon }  \supp Z_a }  \; , 
\end{equation}
where for any subset $B$ in $\bR^d$, $\overline{B}$ stands for its closure.  We also introduce the total occupation measure of $Z$ by setting 
\begin{equation}
\label{defboldM}
 \bM= \int_0^\infty Z_a \; da 
 \end{equation}
whose support is in ${\bf R}$. Next, for any $r \in (0, 1/e)$, we set 
\begin{equation}
\label{defingg}
 g(r)= \frac{r^4}{ (\log \log 1/r)^3 } \; .
 \end{equation}
We denote by $\cP_g$ the $g$-packing measure on $\bR^d$, whose definition is recalled in Section 
\ref{packingsec}. The following theorem is the main result of the paper. 
\begin{theorem}
\label{mainth} Let us assume that $d \geq 5$. Let $g$ be defined by (\ref{defingg}). Fix $\mu \in M_f (\bR^d)$. Then, there exists a positive constant $\kappa_d $ that only depends on space dimension $d$ such that $\bP_\mu$-a.s. for any Borel set $B$ we have 
$$  \bM (B)= \beta \kappa_d  \cdot  \cP_g ( B \cap {\bf R} \, ) \; .$$
\end{theorem}

\vspace{5mm}

We shall actually derive Theorem \ref{mainth} from a similar result concerning the occupation measure of the Brownian Snake that is a process introduced by Le Gall in \cite{LG93} to facilitate the study of super-processes. More precisely, we consider the Brownian Snake $W=(W_t , t \in [0, \sigma])$
with initial value $0$ under its excursion measure denoted by $\bbN_0$. Here, $\sigma$ stands for the total duration of the excursion.  We informally recall that $W$ is a 
continuous Markov process that takes its values in the set of stopped $\bR^d$-valued paths; namely,  under $\bbN_0$ and for any $t \in [0 , \sigma] $, $W_t $ is an application from a random time-interval $[0, H_t]$ to $\bR^d$ such that $W_t (0)= 0$; the process $H= (H_t , t \in [0, \sigma])$ is called the lifetime process of $W$ and it is distributed under $\bbN_0$ as Itô's positive Brownian excursion; 
conditionally given $H$, for any $0 \leq t_1<t_2  \leq \sigma$, the joint law of $(W_{t_1}, W_{t_2})$ is characterised as follows. 
\begin{itemize}
\item $\! W_{t_1}$ is distributed as a $d$-dimensional Brownian path on $[0, H_{t_1}]$ with initial value $0$;
\item $\! W_{t_1} (s)= W_{t_2} (s)$ for any $s \leq m(t_1, t_2):=\inf \{ H_t \, ; \, t_1 \leq t \leq t_2 \}$;

\item $\! (W_{t_2} ( s + m(t_1, t_2)) -W_{t_2} (m(t_1, t_2)) \, ; \, 0 \leq s \leq  H_{t_2} -m(t_1, t_2)\, )$ is distributed as a Brownian path on $[\, 0\, ,\,  H_{t_2} -m(t_1, t_2)\, ]$ with initial value $0$ that is independent from $W_{t_1}$.
\end{itemize}

  For any $t \in [0, \sigma]$, we set $\widehat{W}_t= W_t (H_t)$; $\widehat{W}= (\widehat{W}_t, t\in [0, \sigma])$ is called the endpoint process of $W$. Note that the range of the endpoint process $\widehat{W}$ is a compact subset of $\bR^d$ under $\bbN_0$; we denote it by 
   \begin{equation}
\label{rangeWdef}
 \cR_W = \left\{ \; \widehat{W}_s \; , \; s \in [0, \sigma ] \; \right\} \; .
\end{equation}
The occupation measure of $\widehat{W}$ is the random measure $\cM_W$ given by 
 \begin{equation}
\label{occuWdef}
\langle \cM_W , f \rangle  = \int_0^\sigma f( \widehat{W}_s) \, ds  
\end{equation}
for any positive measurable function $f$ on $\bR^d$. To simplify notation, we simply write $\cR= \cR_W$ and $\cM= \cM_W$ when there is no ambiguity. 
We prove the following results on $\cR$ and $\cM$.  
\begin{theorem} 
\label{typilowdensM} Assume $d \geq 5$. There exists a constant $\kappa_d \in (0, \infty) $ that only depends on space dimension $d$ such that 
\begin{equation}
\label{typidensMM}
\bbN_0 {\rm -a.e.} \;  
{\rm for} \;  \cM \!{\rm -almost} \; {\rm all} \; x \, ,  \, \,  \quad 
 \liminf_{r \rightarrow 0+} \frac{\cM (B(x,r))}{g(r)} = \kappa_d \; .
\end{equation}
\end{theorem}
\begin{theorem}
\label{packingsnake} Assume $d \geq 5$.  $\bbN_0$-a.e. for any Borel set $B$ we have 
$$\cM (B)= \kappa_d \cdot \cP_g (B \cap \cR)  \; .$$ 
\end{theorem}

 Theorem \ref{packingsnake} can be used to get a similar result for Integrated Super-Brownian Excursion measure (ISE). Informally, ISE is a random probability measure $\cM^{(1)}$ that is distributed as $\cM$ under the probability measure $\bbN_0 (\, \cdot \, | \,  \sigma = 1\,   )$. More precisely, we introduce the normalised Brownian Snake $W^{(1)}= (W^{(1)}_t ; t \in [0, 1]) $ as the path-valued process constructed as the Brownian Snake except that its lifetime process $(H^{(1)}_t; t\in [0, 1])$ is distributed as a positive Brownian excursion contidionned to have total duration $1$ (see for instance Bertoin \cite{Be} Chapter VIII-4 for a definition). Then, we set $\widehat{W}^{(1)}_t = W^{(1)}_t (H^{(1)}_t)$ and 
$$ \langle \cM^{(1)} , f \rangle = \int_0^1 f ( \widehat{W}^{(1)}_t ) \, dt \; .$$
We also set $\cR^{(1)}= \{\widehat{W}^{(1)}_t \, ; \, t \in [0, 1]  \}$. Then $\supp \cM^{(1)} \subset \cR^{(1)}$. We easily adapt the proof given by Le Gall \cite{LGRange} p. 313 to derive 
from Theorem \ref{packingsnake} the following result for ISE: if $d \geq 5$, then a.s. for any Borel set $B$ in $\bR^d$ we have 
\begin{equation}
\label{ISEconsq}
\cM^{(1)} (B)= \kappa_d \cdot \cP_g (B \cap \cR^{(1)})  \; .
\end{equation}
(Since the arguments are the same as in \cite{LGRange}, we omit the proof of (\ref{ISEconsq}).)

\medskip

   The paper is organized as follows. In Section \ref{packingsec} we recall the definition of packing measures and useful properties such as the now standard comparison results from Taylor and Tricot \cite{TaTr} (stated as Theorem \ref{TaTrcomparesu}) as well as a more specific density result recalled from Edgar \cite{Edgar07}, that is stated as Lemma \ref{equalitycomp}. In Section \ref{PrelBSnake}, we recall (mostly from Le Gall \cite{LG93}) the definition of Brownian Snake and several path-decompositions that are used in the proof section. In Section \ref{estimatesec},  
we prove several key estimates on the Brownian Snake and the Brownian Tree. Section \ref{proofsec} is devoted to the proof of the  results stated in introduction section: we first prove Theorem \ref{typilowdensM}, then we prove Theorem \ref{packingsnake} from which we derive Theorem \ref{mainth}.

\section{Notation, definitions and preliminary results.}
\label{Prelresu}  

\subsection{Packing measures.}
\label{packingsec}
In this section we gather results concerning packing measures.  We first briefly recall the definition of packing measures on the Euclidian space $\bR^d$. Let $g$ be defined by (\ref{defingg}). 
Let $B$ be any subset of $\bR^d$ and let $\varepsilon \in (0, \infty)$; a closed $\varepsilon$-packing of $B$ is a finite 
collection of pairwise disjoint closed ball $(\overline{B}(x_m, r_m), 1\leq m \leq n)$ whose centers $x_m$ belong to $B$ and whose radii $r_m$ are not greater than $\varepsilon$; we set 
\begin{equation}
\label{premdeltameadef}
 \cP^{(\varepsilon)}_g (B) = \sup \left\{ \sum_{m=1}^n g(r_m) \; ;  \; \left( \overline{B}(x_m, r_m), 1\leq m \leq n \right)   \; \varepsilon \!{\rm -packing} \; {\rm of} \; B \; \right\}  . 
 \end{equation}
and 
\begin{equation}
\label{packingpremeadef}
\cP^*_g (B)= \lim_{\varepsilon \rightarrow 0+}  \cP^{(\varepsilon)}_g (B) \; \in  \; [0, \infty] \; , 
\end{equation}
that is the $g$-packing pre-measure of $B$. The $g$-packing outer measure of $B$ is then given by 
$$ \cP_g (B)= \inf \left\{  \sum_{n \geq 0} \cP^*_g (B_n) \; ; \; B \subset \bigcup_{n \geq 0} B_n \;  \right\} \; .$$
\begin{remark}
\label{slightdif}
The definition of $\cP^{(\varepsilon)}_g$ that we adopt here is slightly different from the definition given by Taylor and Tricot \cite{TaTr} who take the infimum of $\sum_{m=1}^n g(2r_m) $ over $\varepsilon$-packings with open balls. However, since $g$ is a continuous regularly varying function, 
the resulting packing pre-measure $\cP^*_g$ given by (\ref{packingpremeadef}) is $1/16$ times the $g$-packing pre-measure resulting from Taylor and Tricot's definition and the difference is irrelevant for our purpose. \cq  
\end{remark}
We next recall several properties of $\cP_g$ from \cite{TaTr} (see Lemma 5.1 \cite{TaTr}): {\it firstly}, $\cP_g$ is a metric outer measure, all Borel sets are $\cP_g$-measurable and $\cP_g$ is Borel-regular; {\it secondly}, it is obvious from the definition that for any subset $B\subset \bR^d$, we have  
\begin{equation}
\label{premeabound}
\cP_g (B) \leq \cP^*_g (B) \; ; 
\end{equation}
moreover if $B$ is a $\cP_g$-measurable such that $0<\cP_g (B) < \infty$, then for any $\varepsilon >0$, there exists a closed subset $F_\varepsilon \subset B$ such that 
\begin{equation}
\label{innerreg}
\cP_g (B) \leq \cP_g (F_\varepsilon) + \varepsilon \; .
\end{equation}
We also recall here Theorem 5.4 \cite{TaTr} that is a standard comparison result for packing measures.
\begin{theorem} 
\label{TaTrcomparesu} (Theorem 5.4 \cite{TaTr}) Let $\mu$ be a finite Borel measure on $\bR^d$. Let  $B$ be a Borel subset of $\bR^d$. There exists a constant $C>1$ that only depends on space dimension $d$, such that the following holds true. 
\begin{itemize} 
\item{(i)}  If  $ \liminf_{r \rightarrow 0} \frac{\mu (B(x, r))}{g(r )} \leq 1 $ for any $x \in B$, then $  
\cP_g (B)  \geq C^{-1} \mu (B) $. 
\item{(ii)} If  $ \liminf_{r \rightarrow 0} \frac{\mu (B(x, r))}{g(r)} \geq 1 $ for any $x \in B$, then 
$  \cP_g (B)  \leq   C\mu (B) $. 
\end{itemize}
\end{theorem}
We shall actually need the following more specific density results that is due to Edgar (see Corollary 5.10 \cite{Edgar07}). 
\begin{lemma}
\label{equalitycomp} (Corollary 5.10 \cite{Edgar07}) Let $\mu$ be a finite Borel measure on $\bR^d$. Let $\kappa  \in (0, \infty)$ and let $B$ be a Borel subset of $\bR^d$ such that 
$$ \forall x \in B \; , \quad \liminf_{r \rightarrow 0+} \frac{\mu (B(x,r))}{g(r)} = \kappa  \; .$$
Then $\mu (B)= \kappa \cdot \cP_g (B) $.
\end{lemma}
\begin{remark}
\label{StrVitProp} Let us make a brief comment on this result: the main purpose of Edgar's paper \cite{Edgar07} is to deal with fractal measures in metric spaces with respect to possibly irregular gauge functions. Corollary 5.10 \cite{Edgar07} (stated as Lemma \ref{equalitycomp}) holds true in this general setting if $\mu$ satisfies the {\it Strong Vitali Property} (see \cite{Edgar07} p.43 for a definition and a discussion of this topic). Since Besicovitch \cite{Besi45} has proved that any finite measure on $\bR^d$ enjoys the Strong Vitali Property, Lemma \ref{equalitycomp} is an immediate consequence of Edgar's Corollary 5.10 \cite{Edgar07}. \cq 
\end{remark}

\subsection{The Brownian Snake.}
\label{PrelBSnake}
In this section we recall the definition of  the Brownian snake and  the Brownian Tree. We also recall useful properties that are needed in the proof sections. We refer to Le Gall \cite{LG93} or \cite{LG99} for more details. Let us first mention that although we often work on the canonical space for Brownian Snake, we shall sometimes need to introduce an auxiliary measurable space that we denote by $(\Omega, \cG)$ and that is assumed to be sufficiently large to carry a $d$-dimensional Brownian motion denoted by $(\xi_t, t \geq 0; \bP_y, y \in \bR^d)$ as well as the other additional independent random variables we may need. 

\medskip

\noi
$\bullet$ {\it Brownian Snake.} We denote by $\cW$ the set of stopped $\bR^d$-valued paths. A stopped path $w$ in $\cW$ is a continuous application $w: [0, \zeta ] \rightarrow \bR^d$ and the nonnegative number $\zeta= \zeta_w$ is called the lifetime of $w$. The endpoint of $w$ is the terminal value $w(\zeta_\w)$ that is denoted by $\widehat{w}$. We equip $\cW$ with the metric $\delta$ given by 
$$ \delta (w_1, w_2)= \sup_{t \geq 0 }  \;  \lVert  w_1 ( t \wedge \zeta_{w_1} ) - w_2 ( t \wedge \zeta_{w_2})  \rVert  \,   + \,  | \zeta_{w_1}- \zeta_{w_1}|  \; .$$
Then $(\cW, \delta)$ is a separable metric space. Let us fix $x \in \bR^d$. 
We denote by $\cW_x$ the set of stopped paths $w$ such that $w(0)= x$. We identify 
the trivial path $w\in \cW_x$ such that $\zeta_{w} = 0$ with the point $x$ in $\bR^d$.

  The Brownian Snake with initial value $x$ is the strong $\cW_x$-valued continuous Markov process $W=(W_s, s \geq 0)$ that is characterised by the following properties. 
\begin{itemize}
\item{Snake(1):} the lifetime process $\zeta_{W_s}:= H_s$, $s \in [0, \infty)$ is a reflecting Brownian motion. 
\item{Snake(2):} conditionally given the lifetime process  $(H_s , s \geq 0)$, the snake $W$ is distributed as an inhomogeneous 
Markov process whose transitions are described by the following properties: let us fix $s_1< s_2$ and let us set $m (s_1, s_2) := \inf_{u \in [s_1, s_2]} H_u$; then,  
\begin{itemize}
\item{(a)} for any $0 \leq t \leq m (s_1, s_2) $, we have $W_{s_1} (t)= W_{s_2} (t) $; 
\item{(b)} the process $ ( W_{s_2} ( t+ m (s_1, s_2) )  - W_{s_2} (m (s_1, s_2) ) ; 0 \leq t \leq H_{s_2}-m (s_1, s_2) )$ is a standard $d$-dimensional Brownian motion that is independent of $W_{s_1}$.   
\end{itemize}
\end{itemize}

  By convenience, we work on the canonical space of continuous applications from $[0, \infty)$ to $\cW$ that is denoted by $C ([0, \infty) , \cW)$ and $W$ stands for the canonical process.  We denote by $\bbP_x$ the distribution of the Brownian Snake with initial value $x$, and for any $w \in \cW$, we denote by $\bbP_w$ the distribution of the snake with initial value $w$. We also denote by $\bbP^*_w$ the law under 
$\bbP_w$ of $(W_{s\wedge \sigma } , s \geq 0)$ where $\sigma= \inf \{ s >0 \; : \; H_s = 0 \}$.

  Observe that the trivial path $x$ is regular for the Brownian snake. We denote by $\bbN_x$ the excursion measure of $W$ out of state $x$ whose normalisation is specified by: 
\begin{equation}
\label{normexcu}
\bbN_x \left(   \sup_{t \in [0, \sigma ] } H_t > a \right) = \frac{1}{2a} \; , \quad a \in (0, \infty) \; .
\end{equation}

We now recall from \cite{LG93} the connection between Brownian Snake and Super-Brownian motion: recall notation $\cR_W$ and $\cM_W$ from (\ref{rangeWdef}) and (\ref{occuWdef});  let 
$\mu \in M_f (\bR^d)$ and let 
$$ \cQ (dxdW) = \sum_{j \in \cJ} \delta_{ (x_j, W^j) }$$
be a Poisson point process on $\bR^d \times C ([0, \infty) , \cW)$ with intensity $\mu(dx) \bbN_x (dW)$. 
Results due to Le Gall \cite{LG93} entail that there exists a Super-Brownian motion $Z= (Z_t , t \geq 0)$ with branching parameter $\beta= 1$ and initial value $Z_0= \mu$ such that 
\begin{equation}
\label{connection}
{\bf R} \cup \{ x_j \; , \; j \in \cJ \} = \bigcup_{j \in \cJ} \cR_{W^j}  \quad {\rm and} \quad {\bf M} = \sum_{j \in \cJ } \cM_{W^j} \; , 
\end{equation}
where ${\bf R}$ and $ {\bf M}$ are deduced from $Z$ by (\ref{defboldR}) and (\ref{defboldM}) (with an obvious notation for the $\cR_{W^j}$'s and the $\cM_{W^j}$'s). This implies that for any $x \in \bR^d$ and for any nonnegative Borel function $f$:
$$ \bbN_x \left( 1-e^{-\langle \cM , f \rangle  }\right) = -\log \left( \bE_{\delta_x} \left[ \exp (-\langle {\bf M} , f \rangle )\right]\right) \; .$$
Recall that $(\xi_t, t \geq 0; \bP_y, y \in \bR^d)$ stands for a $d$-dimensional Brownian motion defined on the auxiliary measurable space $(\Omega, \cG)$. 
If we denote by $u_f (x)= \bbN_x \left( 1-e^{-\langle \cM , f \rangle  }\right)$, standard results on Super-Brownian motion entail
\begin{equation}
\label{occupequa}
u_f (x) + 2\int_0^\infty  dt \, \bE_x \left[ u_f (\xi_t )^2 \right]  =\int_0^\infty  dt \,  \bE_x \left[ f( \xi_t)  \right] \; .
\end{equation}
(we refer to \cite{LG99} for a proof). Then, an easy argument  implies 
\begin{equation}
\label{occupequa}
\bbN_x \left( \int_0^\sigma ds \, f (\widehat{W}_s ) \right)= \bbN_x \left( \langle \cM , f \rangle \right) = \int_0^\infty  dt \,  \bE_x \left[ f( \xi_t)  \right] \; .
\end{equation}

\medskip

 \noi
$\bullet$ {\it Brownian Tree.} The lifetime process $H=(H_s, 0\leq s \leq \sigma )$ under $\bbN_x$ is distributed as the excursion of the reflecting Brownian motion in $[0, \infty)$. Namely the "law" of $H$ under $\bbN_x$ is It\^o's positive excursion measure of the Brownian motion whose normalisation is given by (\ref{normexcu}); we denote It\^o's positive excursion measure by $N$ and 
we slightly abuse notation by keeping denoting the canonical excursion process under $N$ by $H$.

  The endpoint process of Brownian Snake $\widehat{W}= (\widehat{W}_s, 0 \leq s \leq \sigma)$ can be viewed as a specific coding of the spatial positions of a population combining a branching phenomenon with spatial motion; the lifetime process $H$ is then the contour process of the genealogical tree of the population; this tree is actually  distributed as the Brownian Tree, whose definition in \cite{LG2} (or in \cite{Al2}, in a slightly different context) is given as follows:  for any $s,t \in [0, \sigma] $, we set 
\begin{equation}
\label{distH}
m(s,t) = \inf_{u\in [s\wedge t , s\vee t]} H_u \quad {\rm and} \quad d_H(s,t)= H_t +H_s -2m(s,t)  \; . 
\end{equation}
The quantity $d_H(s,t)$ represents the distance between the points corresponding to $s$ and $t$ in the Brownian Tree. Therefore, two real numbers $t,s \in [0, \sigma] $ correspond to the same point in the Brownian Tree iff $d_H(s, t) = 0$, which is denoted by $s\sim_H t$. Observe that $\sim_H$ is an equivalence relation. The Brownian Tree is given by the quotient set $\cT = [0, \sigma ] / \sim_H$; $d_H$ induces a true (quotient) metric on $\cT$ that we keep denoting $d_H$ and $(\cT, d_H)$ is a random compact metric space that is taken as the definition of 
the Brownian Tree (more specifically,  it is a $\bR$-tree: see \cite{DuLG2} for more details).

  The end point process $\widehat{W}$ can be viewed as a Gaussian process indexed by the Brownian Tree. More precisely,  we recall that there exists a regular version of the conditional distribution of $W$ under 
$\bbN_x$ given the lifetime process $H$. This regular version is a random probability measure on $C([0, \infty), \cW)$ denoted by $Q^H_x $ and we have: 
$$ \bbN_x (dW) = \int N(dH) \, Q^H_x ( dW) \; .$$ 
In view of Property Snake(2),  $\widehat{W}= (\widehat{W}_s, 0 \leq s \leq \sigma)$ under $Q^H_x$ is distributed as a Gaussian process whose covariance is characterized by the following:  
\begin{equation}
\label{covar}
 Q^H_x (\widehat{W_0}= x)= 1 \quad  Q^H_x \left( \lVert \widehat{W_t}-\widehat{W_s} \rVert^2 \right) = 
d_H (s, t) \; , \quad s,t \in [0, \sigma ] \; . 
 \end{equation}
We refer to \cite{DuLG2} for a more intrinsic point of view on spatial trees, namely, $\bR$-trees embedded in $\bR^d$. 

\medskip 

\noi
$\bullet$ {\it Markov property and path-decompositions of $W$.} Markov property for $W$ also applies under $\bbN_x$ as follows: denote by $(\cF_t, t \geq 0)$ the canonical filtration on $C([0, \infty) , \cW)$. Let $T$ be a $(\cF_t, t \geq 0)$-stopping time. Then the law 
of $(W_{T+s} , s \geq 0)$ is $\bbP^*_{W_T} $. Namely for any $\Lambda \in \cF_{T+}$ and for any nonnegative measurable functional $F$, we have 
\begin{equation}
\label{MarkovN}
\bbN_x \left( \un_{ \{ T < \sigma \} \cap \Lambda} F ( W_{T+ s} , 0 \leq s \leq \sigma -T) \right) = \bbN_x \left( \un_{ \{ T < \sigma \}\cap \Lambda} \bbE^*_{W_T} [F]  \right) \; .
\end{equation}
We refer to Le Gall \cite{LG94bis} for more details. 

We shall use (\ref{MarkovN}) in combination with the following 
Poissonnian decomposition: let us fix $w\in \cW_x$; recall notation 
$m (s,t)= \inf \{ H_u \; ; \;  s\wedge t \leq u  \leq s \vee t  \}$. To avoid trivialities, we assume that $\zeta_w >0$. Observe that $\bbP^*_w$-almost surely, for any $s \in [0, \sigma)$, $W_s (t)= w(t)= W_0 (t)$ for any  
$t\in [0,  m(0, s)]$. We keep using notation $H_s= \zeta_{W_s}$ for the lifetime process. Let us denote by 
$(l_i, r_i)$, $i \in \cJ$ the excursion intervals of the process $(H_s-m(0,s), s \in [0, \sigma ])$ above $0$. For any $i \in \cJ$,  and for any $s \geq 0$, we set 
$$ H^i_s = H_{ (l_i +s )\wedge r_i } -H_{l_i} \quad {\rm and} \quad W^i_s (t) = W_{(l_i +s )\wedge r_i }  (H_{l_i} +t ) \; , \; t \in [0, H^i_s]  \; .$$  
Then, we recall from Le Gall  \cite{LG94bis} the following property: under $\bbP^*_w$, the point measure 
\begin{equation}
\label{Poisdecinit} 
\cN(dt dW)  = \sum_{i \in \cJ } \delta_{(H_{l_i} , W^i  )}
\end{equation}
is a Poisson point measure with intensity $2\cdot \un_{[0, \zeta_w]} (t) dt \bbN_{w(t)} (dW) $. This decomposition combined with Markov property under $\bbN_x$ implies that for any $(\cF_t, t \geq 0)$-stopping time $T$, any nonnegative Borel measurable function $f$ and any nonnegative measurable functional $F$,  we have 
 \begin{eqnarray}
\label{PoisdecMarkov} 
\bbN_x \left( \un_{ \{ T < \sigma \}} F( W_{\cdot \wedge T }) \exp \left( - \int_T^\sigma ds f(  \widehat{W}_s  ) \right)  \right) = 
\hspace{50mm} \nonumber \\
\bbN_x \left( \un_{ \{ T < \sigma \}}  F( W_{\cdot \wedge T }) 
\exp \left(   - 2   \int_0^{H_T} dt \; \bbN_{_{W_{_T} (t)}} \! \!  \left( 1-e^{-  \langle \cM , f \rangle} \right) \right)   \right) \; .
\end{eqnarray}

  We shall apply (\ref{PoisdecMarkov}) at deterministic times and at hitting times of closed balls that are discussed here: for any $x \in \bR^d$ and for any $r >0$, we denote by $B(x, r)$ the open ball with center $x$ and radius $r$ and we write $\overline{B} (x, r)$ for the corresponding closed ball. For any $w \in \cW$ we set 
\begin{equation}
\label{notahitting}
\tau_{x, r} (w)= \inf \{  t \in [0, \zeta_w ]   \; : \;   w(t) \in \overline{B} (x, r)  \} \; , 
\end{equation}
with the convention: $\inf \emptyset = \infty$. Then, $  \tau_{x, r} ( \widehat{W} ) $ is a $(\cF_t ,t \geq 0)$-stopping-time and observe that  $\tau_{x, r} ( \widehat{W} ) < \infty$ iff $\cR \cap \overline{B} (x, r) \neq \emptyset$. 

\medskip

  We next set: 
\begin{equation}
\label{uxrdef}
 u_{x,  r} (y) = \bbN_y  \left( \cR \cap \overline{B} (x, r) \neq \emptyset \right) \; .
 \end{equation}
We need to recall several important properties of $u_{x,r}$ whose proofs can be found in Le Gall \cite{LG99}: first of all,  $u_{x,  r}$ is twice continuously differentiable in $\bR^d \backslash  \overline{B} (x, r) $ and it satisfies 
\begin{equation}
\label{Dirichlet}
 \Delta u_{x , r} (y) = 4  u_{x , r}^2 (y) \; , \quad  y \in \bR^d \backslash  \overline{B} (x, r) \; .
\end{equation}
Next, $u_{x, r} (y) \rightarrow \infty$ when $\lVert y -x \rVert$ goes to $r$; since $\cR$ is compact, we also have $u_{x,r} (y) \rightarrow 0$ when $\lVert y \rVert$ goes to $\infty$. Moreover $u_{x,r}$ is the maximal nonnegative solution of (\ref{Dirichlet}).

  Let us briefly discuss further (simple) properties of $u_{x,r}$ that are needed in the proofs section: a symmetry argument first implies that 
\begin{equation}
\label{Symmetry}
u_{x , r} (y) =  u_{0 , r} (y-x) \; , \quad  y \in \bR^d \backslash  \overline{B} (x, r) \; .
\end{equation}
Next, observe that $u_{0 , r}$ is radial. Namely, there exists a twice continuously differentiable application $v_r : (r, \infty) \rightarrow (0, \infty) $ such that $u_{0, r} (y)= v_r ( \lVert y \rVert ) $; moreover by (\ref{Dirichlet}) $v_r$ is the unique solution of the following ordinary differential equation:
\begin{equation}
\label{radial}
 v''_r (t) + \frac{d-1}{t} v'_r (t) = 4 v_r^2 (t) \; , \; t \in (r, \infty)  \quad  {\rm with} \quad  \lim_{t \downarrow r} v_r (t) = \infty \; , \; \lim_{t \uparrow \infty} v_r (t) = 0 \; . 
\end{equation}
To simplify notation, we set $u:= u_{0, 1}$ and $v:=v_{0, 1}$. Namely, 
\begin{equation}
\label{notau1v1}
v (\lVert y \rVert ) = u(y)= u_{0, 1} (y)= \bbN_y   \left( \cR \cap \overline{B} (0, 1) \neq \emptyset \right) \; , \quad  y \in \bR^d \backslash  \overline{B} (0, 1) \; .
\end{equation}
Ordinary differential equation (\ref{radial}) implies  $v_r (t)= r^{-2} v (r^{-1} t)$. Therefore we get  
\begin{equation}
\label{Scaling}
u_{x , r} (y) =  r^{-2} u (r^{-1}( y-x)) \; , \quad  y \in \bR^d \backslash  \overline{B} (x, r) \; .
\end{equation}
Finally, the maximum principle easily entails that $v(t) \leq v(2)2^{d-2} t^{2-d} $, for any $t \geq 2$. Thus, by (\ref{Scaling}) we get:  
\begin{equation}
\label{boundu}
u_{x,  r} (y) \leq  v(2)2^{d-2} r^{d-4} \lVert y-x \rVert^{2-d} \; , \quad y \in \bR^d \backslash \overline{B} (x, 2 r)\; .
\end{equation}
This upper bound shall be often used in the proofs.

\vspace{5mm}

   We next describe the distribution of the snake when it hits for the first time a closed ball. Recall notation $(\Omega, \cG)$ and 
$(\xi_t, t \geq 0; \bP_y, y \in \bR^d)$. Let us fix $x,y\in \bR^d$ and $R>r>0$ such that $\lVert y-x\rVert >R$. To simplify notation, we set 
$$ \tau_{x,r} (\widehat{W} ) = \tau \; .$$ 
A result due to Le Gall \cite{LG94} (see also \cite{DuLG}, Chapter 4) asserts that  $W_{\tau}$ under the probability measure $\bbN_y \left( \, \cdot \, | \, \cR \cap \overline{B} (x,r)\neq \emptyset \right)$ is distributed as the solution of the stochastic equation 
$$ dX_t = d\xi_t + \frac{\nabla u_{x,r}}{u_{x,r}} (X_t) dt \; , \quad X_0= y \; , $$
where we recall that $\xi$ stands for a $d$-dimensional Brownian motion. By applying Girsanov's theorem, we can prove that  for any nonnegative measurable functional 
$F$ on $\cW$, we have 
\begin{eqnarray}
\label{charachitti}
\bbN_y \left( \un_{ \{  \cR \cap \overline{B} (x,r)\neq \emptyset   \}  } F \left(   \, W_{\tau  } ( t ) \, ;  \, 0 \leq t \leq \tau_{x, R} (W_\tau )  \,\right) \;  \right) =\hspace{50mm} \nonumber \\ 
\bE_y \left[ u_{x, r} (\xi_{\tau_{x, R} (\xi) }) \, F( \xi_t ; 0 \leq t \leq \tau_{x, R} (\xi) \, ) \, . \, e^{-2\int_0^{\tau_{x, R} (\xi) }  u_{x, r} (\xi_s) } \right] \; . 
\end{eqnarray}
We refer to \cite{DuLG}, Chapter 4, pp.131-132 for a proof of this specific result which is only used in Lemma \ref{endestimate}.

\medskip

\noi
$\bullet$ {\it Palm decomposition of Brownian Snake occupation measure.} The proof of Theorem \ref{typilowdensM} heavily relies on the following Palm decomposition of $\cM$ whose proof can be found in Le Gall \cite{LG94bis}: recall that $(\xi_t, t \geq 0; \bP_y , y \in \bR^d)$ stands for a $d$-dimensional Brownian motion defined on the auxiliary measurable space $(\Omega, \cG)$. To simplify notation, we assume that it is possible to define on $(\Omega, \cG)$ a point measure on $[0, \infty) \times C([0, \infty), \cW)$ denoted by
\begin{equation}
\label{notapalm}
 \cN^* ( dtdW) = \sum_{j\in \cJ^*} \delta_{ (t_j, W^j)} 
 \end{equation}
whose distribution conditionally given $\xi$ under $\bP_0$ is the distribution of a Poisson point measure with intensity $4\, dt \, \bbN_{\xi(t)} (dW)$. For any $j \in \cJ^*$, we denote by $\cM_j$ the occupation measure of the endpoint process $\widehat{W}^j$ and for any $ a\in (0, \infty)$ we set 
\begin{equation} 
\label{Mstardef}
\cM^*_a = \sum_{j \in \cJ^*} \un_{ [0, a] } (t_j) \cM_j \; .
\end{equation}
Then, for any $x \in \bR^d$ and for any nonnegative measure functional $F$ we have 
\begin{equation}
\label{Palm}
 \bbN_x  \left( \int \! \!  \cM (dy ) \,  F \left( \cM ( B(y,r) ) \, ; \, r \geq 0 \right) \,  \right) = \int_0^{\infty} \! \! \! da \, \bE_0 
 \left[ \, F \left( \cM^*_a (B(0, r) ) \, ; \, r \geq 0 \right) \, \right] . 
\end{equation}

\subsection{Estimates.}
\label{estimatesec}
In this section we prove key estimates used in the proof sections. We first state a result concerning the Brownian Tree: recall that $N$ stands for Itô's excursion measure of Brownian motion and recall that $H= (H_t, 0, \leq t \leq \sigma )$ denote the generic excursion. We assume that the normalisation of $N$ is given by (\ref{normexcu}). Recall notation $d_H$ from (\ref{distH}). We denote the Lebesgue measure on the real line by $\ell$. For any $r \in (0, 1/e)$, we set 
\begin{equation}
\label{defkk}
k(r)= \frac{r^2}{\log \log 1/r}  \; .
\end{equation}
We first prove the following lemma. 
\begin{lemma}
\label{typidensBtree} $N$-almost everywhere, for $\ell$-almost all $t \in [0, \sigma]$, we have  
\begin{equation}
\label{typidenBtreeeq}
\liminf_{r \rightarrow 0+} \;  \frac{1}{k(r)} \int_0^\sigma \un_{ \left\{ d_H (s,t) \leq r  \right\}} \, ds \; \;  \geq  \frac{1}{4}  \; .
\end{equation}
\end{lemma}
\noi
{\bf Proof:} for any $t \in [0, \sigma]$ and any $r >0$, we set
$$ {\bf a}(t,r)=  \int_0^\sigma \un_{ \left\{ d_H (s,t) \leq r  \right\}} \, ds  \; .$$
We prove (\ref{typidenBtreeeq}) thanks to Bismut's decomposition of Brownian excursion: suppose that $(B_t,t\geq 0)$ and $(B'_t,t\geq 0)$ are two $\bR$-valued processes defined on 
$(\Omega, \cG)$ and assume that $\Pi$ is a probability measure under which $B$ and $B'$ are distributed as two independent linear Brownian motions with initial value $0$; 
for any $a >0$, we set
$$T_a=\inf\{t\geq 0: B_t=-a\}\ ,\ T_a'=\inf\{t\geq 0: B'_t=-a\} \; .$$
Then for any nonnegative measurable functional $F$
on $C([0, \infty) ,\bR)^2$, we have
\begin{equation}
\label{Bismutun}
N \left(  \int_0^\sigma  F\left( H_{(t-\cdot)_+} ; \, H_{(t+\cdot)\wedge \sigma} \right) \right)
= \int_0^\infty da\,\Pi \left[ F \left( a+ B_{\cdot \wedge T_a}  \, ; \, a+B'_{\cdot \wedge T'_a}  \right) \right].
\end{equation}
This identity is known as Bismut's decomposition of the Brownian excursion (see \cite{Bis85} or see Lemma 1 \cite{LG1} for a simple proof). For any $r,a >0$, we also set 
$$ {\bf b}(r, a)= \int_0^{T_a} \un_{ \left\{ B_t-2I_t  \leq r  \right\}} \, ds  \quad {\rm and} \quad  {\bf b}'(a,r)= \int_0^{T'_a} \un_{ \left\{ B'_t-2I'_t  \leq r  \right\}} \, ds  \; , $$
where $I_t$ and $I'_t$ stand respectively for $\inf_{s \in [0, t]} B_s$ and $\inf_{s \in [0, t]} B'_s$. Then, (\ref{Bismutun}) implies that for any functional $F$ on $C([0, \infty) ,\bR)$: 
\begin{equation}
\label{Bismutdeux}
N \left(  \int_0^\sigma dt \,  F\left( {\bf a}(t,r) \, ; \, r>0 \right) \right)
= \int_0^\infty da\; \Pi \left[ \,  F \left( \, {\bf b}(r,a) + {\bf b}'(r,a) \, ; \, r >0 \, \right) \, \right].
\end{equation}
Now observe that if $r <a$, then ${\bf b}(r,a)$ and ${\bf b}'(r,a)$ do not depend on $a$. So we simply denote them by ${\bf b}(r)$ and ${\bf b}'(r)$. Therefore, we only need to prove that 
\begin{equation}
\label{Bismuttrois}
 \Pi\! {\rm -a.s.} \quad \liminf_{r\rightarrow 0} \frac{{\bf b}(r)+{\bf b}'(r)}{k(r)} \geq \frac{1}{4} \; .
 \end{equation}
By a famous result due to Pitman, the process $(B_t-2I_t, t\geq0)$ is distributed under $\Pi$ as the 
three-dimensional Bessel process. Therefore, if $a>r$, then 
${\bf b}(r)$ is distributed as the three-dimensional Brownian occupation measure of the unit ball with center $0$. A result due to Ciesielski and Taylor \cite{CiTa} asserts that ${\bf b}(r)$ is distributed as the first exit time of $B$ from interval $[-r,r]$:
$$ \theta_r = \inf \{ t \geq 0 \; : \; \vert B_t \rvert = r \} \; $$
(see also \cite{ReYo} Chapter XI). A standard martingale argument allows to explicitly compute the Laplace transform of $\theta_r$: for any $\lambda >0$, we have 
\begin{equation}
\label{exitinterv}
\Pi \left[ \exp (-\lambda \theta_r) \right]= ( \cosh (r\sqrt{2\lambda}))^{-1}.
\end{equation}
Thus, we get 
$$ \Pi \left[ e^{- \lambda ({\bf b}(r)+{\bf b}'(r)) } \right] =(\cosh (r\sqrt{2\lambda}))^{-2}  \leq 4e^{-r\sqrt{8\lambda}} \; .$$
Therefore, for any $\lambda >0$, Markov inequality entails 
$$ \Pi \left[ \,   {\bf b}(r)+{\bf b}'(r) \leq k(r)  \, \right] \leq 4 \exp ( \varphi (\lambda ) ) \; , $$
where we have set $\varphi (\lambda )= \lambda (\log \log 1/r)^{-1} - \sqrt{8 \lambda}$; $\varphi$ reaches its minimal value at $\lambda_0 = 2 (\log \log 1/r)^2$ and $\varphi(\lambda_0) = -2 \log \log 1/r$. Consequently, 
$$ \sum_{n \geq 0}  \Pi \left[ \,   {\bf b}(2^{-n})+{\bf b}'(2^{-n}) \leq k(2^{-n}) \, \right] \; < \; \infty \; , $$
which easily entails (\ref{Bismuttrois}) and which completes the proof of the lemma. \cqfd



\medskip

We next provide an estimate on the tail distribution at $0+$ of $\cM (B(0, r))$ under $\bbN_0$ (Lemma \ref{endestimate}). To that end we state the following preparatory lemma that is only used in the proof of Lemma \ref{endestimate}. 
\begin{lemma} 
\label{loweroccup}
There exists a constant $C_1  \in (0, \infty)$ that only depends on space dimension $d$, such that for any $r >0$, and for any $\lambda >1$, the following inequality holds true. 
$$ \bbN_0 \left(  1-e^{-r^{-4}\lambda \, \cM (B(0, r)\, )} \right) \geq C_1 r^{-2}\sqrt{\lambda} \; .$$
\end{lemma}
\noi
{\bf Proof:} to simplify notation we set $q(\lambda, r)= \bbN_0 \left(  1-e^{-\lambda \, \cM (B(0, r)\, )} \right)$.   First observe that (\ref{PoisdecMarkov}), combined with elementary arguments, implies 
\begin{eqnarray*}
q(\lambda, r) & =& \bbN_0 \left(
  1-e^{ -\lambda \int_0^\sigma ds \, \un_{ \{ \lVert \widehat{W}_s \rVert \leq r  \} } } \;  \right) \\
 & =&\lambda  \int_0^\infty  dt \,  \bbN_0 \left(   \un_{ \{ t \leq \sigma \; ; \;    \lVert \widehat{W}_t \rVert \leq r  \} }e^{ -\lambda \int_t^\sigma    ds \, \un_{ \{ \lVert  \widehat{W}_s \rVert \leq r  \} }   }\right) \\
&=& \lambda   \int_0^\infty  dt \, \bbN_0  \left(   \un_{ \{ t \leq \sigma \; ; \;    \lVert \widehat{W}_t \rVert \leq r  \} } e^{ - \int_0^{H_t}   ds \,  \bbN_{W_t (s) } \left( 1-e^{ -\lambda \cM (B(0, r))  }  \right)   }\; \right) .
\end{eqnarray*}
Now recall that $ \cM (B(0, r)) \leq \langle \cM , \un \rangle = \sigma $. Thus, for any $y \in \bR^d$, we have 
$$ \bbN_y \left( 1-e^{ -\lambda \cM (B(0, r))} \right)  \leq N \left( 1-e^{-\lambda \sigma  }\right)  .$$
Then, a standard argument in fluctuation theory  asserts that $N \left( 1-e^{-\lambda \sigma  }\right) =
 \sqrt{\lambda /2}  $. Thus, the latter inequality, combined with (\ref{occupequa}) implies:  
\begin{eqnarray*}
 q(\lambda, r) & \geq & \lambda   \int_0^\infty  dt \, \bbN_0  \left(   \un_{ \{ t \leq \sigma \; ; \;    \lVert \widehat{W}_t \rVert \leq r  \} }\exp \left( - H_t  \sqrt{\lambda /2}    \right) \right),  \\
& \geq & \lambda   \int_0^\infty  dt \, \bP_0 ( \lVert \xi_t \rVert  \leq r  )  \, e^{-t \sqrt{\lambda /2}  } , \\
& \geq &  \sqrt{2\lambda}  \int_0^\infty  dt \, \bP_0 ( \lVert \xi_t \rVert   \leq 2^{-1/4} \lambda^{1/4} r ) \,  e^{-t} .
\end{eqnarray*}
By replacing $\lambda $ by $r^{-4}\lambda$ in the previous inequality, the desired result holds true  with $C_1= \sqrt{2}\int_0^\infty dt \, e^{-t}\bP_0 ( \lVert \xi_t \rVert   \leq 2^{-1/4}  )$. \cqfd 

\medskip 

Recall notation $g$ from (\ref{defingg}) and notation $v$ from (\ref{notau1v1}). From Lemma \ref{loweroccup}, we derive the following estimate that is used in the proof of Theorem \ref{packingsnake}. 
\begin{lemma} 
\label{endestimate}
There exist $C_2, C_3, \kappa_0, r_0 \in (0, \infty)$ that only depend on space dimension $d$, and such that for any $r \in (0, r_0) $, for any $\kappa \in (0, \kappa_0 )$, and for any $x \in \bR^d \backslash \overline{B}(0, 2r)$, the following inequality holds true. 
$$ \bbN_0 \left( \;  \cR \cap \overline{B} (x,r) \neq \emptyset  \; ; \; \cM (B(x, 2r)) \leq \kappa g(r) \;  \right)  \; \leq C_2 r^{d-4} \lVert  x\rVert^{2-d}  \left( \log 1/r \right)^{- (\kappa/C_3) ^{-1/3}}  \; .$$
\end{lemma}
\noi
{\bf Proof:}  recall notation $\tau_{x,r}$ from (\ref{notahitting}). To simplify notation we set $\tau= \tau_{x,r} (\widehat{W})$. First note that on the event $ \cR \cap \overline{B} (x,r) \neq \emptyset $, we have 
$$  \int_{\tau}^\sigma  \un_{ \{  \lVert \widehat{W}_s -x \rVert < 2r   \}} ds \; \leq \;  \cM (B(x, 2r)) \; .$$
Let us fix $\lambda >0 $. The previous inequality combined with (\ref{PoisdecMarkov}) implies 
\begin{eqnarray}
\label{Laplaceun}
\bbN_0 \left(  \un_{ \{  \cR \cap \overline{B} (x,r) \neq \emptyset  \} } e^{-\lambda  \cM (B(x, 2r))} \right) \leq \hspace{60mm} & &  \nonumber \\ 
\bbN_0 \left(  \un_{ \{  \cR \cap \overline{B} (x,r) \neq \emptyset  \} }
\exp \left(   - 2   \int_0^{H_{\tau}  } dt \; \bbN_{_{W_{_\tau} (t)}} \! \!  \left( 1-e^{  -\lambda \cM (B(x,2r)) }  \right) \right)   \right) \; .
\end{eqnarray}
We want to modify the second member in (\ref{Laplaceun}) in order to get an upper bound thanks to (\ref{charachitti}). To that end, we introduce two times $T_1$ and $T_2$ that are defined as follows: for any $w \in \cW_0$, we set $T_1 (w)=  \tau_{x, \frac{3}{2}r} (w) $. If $T_1 (w)= \infty $, then we set 
 $T_2 (w)= \infty$; if  $T_1 (w)<  \infty $, then we set 
\begin{equation}
\label{T2def}
 T_2 (w)= \inf \left\{ s \in [\, 0\, ,  \; \zeta_w -T_1 (w) \, ] \; : \; \lVert w( s + T_1 (w)) -w(T_1(w))  \rVert   > r/4 \; \right\},  
 \end{equation}
with the usual convention: $\inf \emptyset = \infty$. To simplify notation we also set 
$$T_3 (w)= \tau_{x, \frac{5}{4}r} (w)  \; .$$ 
Observe that $T_1 (w)+T_2 (w) \leq T_3 (w)$. Since $x \in \bR^d \backslash \overline{B}(0, 2r)$, $\bbN_0$-a.e. on the event 
$ \{ \cR \cap \overline{B} (x,r) \neq \emptyset \}$, we have $T_1 (W_\tau)+T_2 (W_\tau) \leq T_3 (W_\tau) < \tau <  \infty $. Moreover, for any 
$ t \in [T_1 (W_\tau) , T_1 (W_\tau)+T_2 (W_\tau)  ]$, the following inequality holds true: 
$$ \bbN_{_{W_{_\tau} (t)}} \! \!  \left( 1-e^{  -\lambda \cM (B(x,2r)) }\right)  \geq \bbN_{_{W_{_\tau} (t)}} \! \!  \left( 1-e^{  -\lambda \cM (B(W_\tau(t) ,r/4)) }\right)  = \bbN_0 \left( 
1-e^{  -\lambda \cM (B(0,r/4)) } \right) . $$ 
Inequality (\ref{Laplaceun}) and Lemma \ref{loweroccup} then entail for any $\lambda >1$, 
\begin{equation}
\label{Laplacedeux}
\bbN_0 \left(  \un_{ \{  \cR \cap \overline{B} (x,r) \neq \emptyset  \} } e^{-r^{-4}\lambda  \cM (B(x, 2r))} \right) \leq \bbN_0 \left(  \un_{ \{  \cR \cap \overline{B} (x,r) \neq \emptyset  \} } e^{ -2C_1r^{-2} \sqrt{\lambda} T_2 (W_\tau)} \right) .
\end{equation}
Then, we set $\mu= 2C_1r^{-2} \sqrt{\lambda} $ and $a=  \bbN_0 \left(  \un_{ \{  \cR \cap \overline{B} (x,r) \neq \emptyset  \} } e^{ -\mu T_2 (W_\tau)} \right)$. We apply (\ref{charachitti}) with $y=0$, $R= 5r/4$ and 
$F(w)= \exp (-\mu T_2 (w))$, and we get the following: 
\begin{eqnarray*}
 a &= & \bE_0 \left[  \un_{ \{ T_3 (\xi) < \infty \}   } u_{x,r} (\xi_{T_3 (\xi)}) \,  e^{-\mu T_2 (\xi) -2\int_0^{T_3 (\xi) }  u_{x,r} (\xi_s) ds }  \right] \\
 & \leq & \bE_0 \left[  \un_{ \{ T_1 (\xi) < \infty \}   } u_{x,r} (\xi_{T_1 (\xi)}) e^{ -2\int_0^{T_1 (\xi) }  u_{x,r} (\xi_s) ds} \cdot 
  \un_{ \{ T_3 (\xi) < \infty \}   } \frac{ u_{x,r} (\xi_{T_3 (\xi)}) }{ u_{x,r} (\xi_{T_1 (\xi)})  } e^{-\mu T_2 (\xi) }  \right] \; .
 \end{eqnarray*}
Recall notation $v$ from (\ref{notau1v1}).  Now observe that (\ref{Symmetry}) and (\ref{Scaling}) imply 
$$\bP_0 {\rm -a.s.} \; {\rm  on} \;  \{ T_3 (\xi) < \infty \}\, ,  \quad 
 \frac{ u_{x,r} (\xi_{T_3 (\xi)}) }{ u_{x,r} (\xi_{T_1 (\xi)})  }= \frac{v_{0, r} (5r/4)}{v_{0, r} (3r/2)}= \frac{v(5/4)}{v(3/2)} := C_4. $$
Note that $C_4$ only depends on space dimension $d$. Then, we get 
\begin{equation}
\label{Laplacetrois}
 a \; \leq  \;   C_4 \,  \bE_0 \left[  \un_{ \{ T_1 (\xi) < \infty \}   } u_{x,r} (\xi_{T_1 (\xi)}) e^{ -2\int_0^{T_1 (\xi) }  u_{x,r} (\xi_s) ds} \cdot 
 e^{-\mu T_2 (\xi) }  \right] \; . 
 \end{equation}
We next apply Markov property at $T_1 (\xi)$ in the right member of the previous inequality: the very definition (\ref{T2def}) of $T_2 (\xi)$, combined with an elementary argument entails the following. 
\begin{eqnarray*}
\bE_0 \left[  \un_{ \{ T_1 (\xi) < \infty \}   } u_{x,r} (\xi_{T_1 (\xi)}) e^{ -2\int_0^{T_1 (\xi) }  u_{x,r} (\xi_s) ds} \cdot 
 e^{-\mu T_2 (\xi) }  \right]  =  \hspace{50mm} & &  \hspace{50mm} \\ 
 \bE_0 \left[  \un_{ \{ T_1 (\xi) < \infty \}   } u_{x,r} (\xi_{T_1 (\xi)}) e^{ -2\int_0^{T_1 (\xi) }  u_{x,r} (\xi_s) ds} \right]  \, \bE_0 \left[  
 e^{-\mu \tau_{0, r/4} (\xi) }  \right] \; . 
 \end{eqnarray*}
If we now apply (\ref{charachitti}) with $y=0$, $R= 3r/2$ and $F= \un$, then we get 
$$  \bE_0 \left[  \un_{ \{ T_1 (\xi) < \infty \}   } u_{x,r} (\xi_{T_1 (\xi)}) e^{ -2\int_0^{T_1 (\xi) }  u_{x,r} (\xi_s) ds} \right]= \bbN_0 \left(  \cR \cap \overline{B} (x,r) \neq \emptyset  \right) = u_{x,r} (0) \; .$$
Then (\ref{boundu}) and (\ref{Laplacetrois}) imply that for any $\lambda >1$, any $r >0$, and any $x \in \bR^d \backslash \overline{B} (0, 2r) $, we have: 
\begin{equation}
\label{Laplacequatre}
a \; \leq \; C_5 r^{d-4} \lVert  x\rVert^{2-d}  \bE_0 \left[  
 e^{-\mu \tau_{0, r/4} (\xi) }  \right]\; , 
\end{equation}
with $C_5:= C_4 v(2)2^{d-2}$, that only depends on space dimension $d$. 

\medskip

We now provide an upper bound for  $\bE_0 \left[  e^{-\mu \tau_{0, r/4} (\xi) }  \right] $: for any $j \in \{ 1, \ldots , d \}$, we denote by $ \xi^{(j)} = (\xi^{(j )}_t , t \geq 0)$ the $j$-th component of $\xi$ 
in the canonical basis of $\bR^d$. Then, observe that $\bP_0$-a.s. 
$$ \tau_{0, r/4}  (\xi)  \geq \min_{1 \leq j \leq d } \inf \{ t \geq 0 \; : \;  \sqrt{d} \, \lvert  \xi^{(j) }_t \rvert  > r/4 \;  \} \; .$$
An easy argument based on (\ref{exitinterv}) implies 
$$ \bE_0 \left[ e^{-\mu \tau_{0, r/4} (\xi) }\right] \leq 2d \cdot e^{-\frac{r}{\sqrt{8d}}  \sqrt{\mu}} \; .$$
We now set $C_2 = 2dC_5$ and $C_6= \sqrt{C_1/ 4d} $. Then (\ref{Laplacedeux}) and  (\ref{Laplacequatre}) imply that for any $\lambda >1$, any $r >0$ and any $x \in \bR^d \backslash \overline{B} (0, 2r) $, 
\begin{equation}
\label{Laplacecinq}
\bbN_0 \left(  \un_{ \{  \cR \cap \overline{B} (x,r) \neq \emptyset  \} } e^{-r^{-4}\lambda  \cM (B(x, 2r))} \right) \; \leq \; C_2 \, r^{d-4} \lVert  x\rVert^{2-d}  
\cdot  e^{- C_6  \, \lambda^{1/4} }  \; .
\end{equation}
To simplify notation we set $\phi (r) = \log \log 1/r $. Thus $g(r)= r^4 \phi (r)^{-3}$. We also set 
$$ b(x, r, \kappa)= \bbN_0 \left( \;   \cR \cap \overline{B} (x,r) \neq \emptyset  \; ; \; \cM (B(x, 2r)) \leq \kappa g(r) \;  \right) \; .$$
Markov inequality and  (\ref{Laplacecinq}) imply that the following inequalitites hold true for any $\lambda >1$, any $r>0$ and any $x \in \bR^d \backslash \overline{B} (0, 2r) $: 
\begin{eqnarray}
\label{finalboundfinal}
 b(x, r, \kappa) & = & \bbN_0 \left( \;   \cR \cap \overline{B} (x,r) \neq \emptyset  \; ; \; r^{-4} \lambda \cM (B(x, 2r)) \leq \kappa \lambda \phi (r)^{-3}  \; \right) \nonumber  \\
 & \leq & e^{\kappa \lambda \phi (r)^{-3} } \bbN_0 \left(  \un_{ \{ \cR \cap \overline{B} (x,r) \neq \emptyset   \} }  e^{-r^{-4}\lambda  \cM (B(x, 2r))} \right)  \nonumber  \\
 & \leq &  C_2 \, r^{d-4} \lVert  x\rVert^{2-d}  
\cdot  e^{ \psi (\lambda) }   \; ,
\end{eqnarray}
where, $\psi (\lambda)=  \kappa \lambda \phi (r)^{-3} -  C_6  \, \lambda^{1/4} $; $\psi $ reaches its minimal value on $[0, \infty)$ at $\lambda_0= (C_6/4\kappa)^{4/3} \phi(r)^4$ and 
$\psi (\lambda_0) = - (\kappa/C_3) ^{-1/3} \phi (r) $, where $C_3= 3^{-3} (4/C_6)^{4}$. If we take $ r_0= e^{-e}$ and $\kappa_0= C_6 /4$, then for any $r <r_0$ and any $\kappa < \kappa_0 $, we have $\lambda_0 >1$ and (\ref{finalboundfinal}) applies with $\lambda= \lambda_0$, which completes the proof of the lemma. \cqfd 

\vspace{5mm}

We shall need the following bound in the proof of Theorem \ref{typilowdensM}. Recall notation 
$u_{0, r}$ and recall that $\xi$ under $\bP_0$ is a standard $d$-dimensional Brownian motion starting at $0$. 
\begin{lemma} 
\label{expectgreen} Assume that $d \geq 5$. Let $b$ and $r $ be positive real numbers such that $b\geq 2r$. 
There exists $C_7\in (0, \infty)$ that only depends on space dimension $d$ such that the following inequality holds true. 
$$ \bE_0 \left[  \int_0^\infty \un_{ \{  \lVert \xi_t \rVert \geq b  \}}  u_{0, r} ( \xi_t) \, dt \;  \right]  \leq  C_7 \left( \frac{r}{b} \right)^{d-4} \; .$$
\end{lemma}
\noi
{\bf Proof:} by (\ref{boundu}), by the scaling property of Brownian motion and thanks to a spherical change of variable, we get the following inequalities. 
\begin{eqnarray}
\label{ineqbess}
 \bE_0 \left[  \int_0^\infty \un_{ \{  \lVert \xi_t \rVert \geq b  \}}  u_{0, r} ( \xi_t) \, dt \;  \right] & \leq & 
 v(2)2^{d-2} r^{d-4} \int_0^\infty \bE_0 \left[  \un_{ \{  \lVert \xi_t \rVert \geq b  \}} \lVert  \xi_t\rVert^{2-d} \right] \, dt,  \nonumber  \\
 & \leq & C_8 r^{d-4} \int_b^\infty d \rho \, \rho \int_0^\infty  dt \,  t^{-d/2} e^{-\frac{\rho^2}{2t}}   \; , 
  \end{eqnarray}
where $C_8 \in (0, \infty)$ only depends on space dimension $d$. We next use the change of variable $s: = \rho^2/(2t)$:  then, there exists 
$C_9 \in (0, \infty)$ that only depends on $d$ such that 
\begin{eqnarray*}
 \int_b^\infty d \rho \,  \rho \int_0^\infty  dt \, t^{-d/2} e^{-\frac{\rho^2}{2t}} & = & C_9 \int_b^{\infty} d \rho \, \rho^{3-d} \int_0^\infty ds \, s^{\frac{d}{2} -2} e^{-s} \\
  &=& \frac{C_9 \Gamma (\frac{d}{2} -1)}{d-4} b^{4-d}  \; .
  \end{eqnarray*}
This inequality combined with (\ref{ineqbess}) entails the desired results with $C_7=   \frac{C_8 C_9 \Gamma (\frac{d}{2} -1)}{d-4}$. \cqfd

\section{Proof of the results. }
\label{proofsec}

\subsection{Proof of Theorem \ref{typilowdensM}.}
\label{prooftypilowdensMsec}
Recall notation $\cN^* (dtdW)$ from (\ref{notapalm}). $\cN^* (dtdW)$ defines a collection of random points $\{ (t_j, W^j) \, ; \, j \in \cJ^*\} $ in $[0, \infty)\times \cW$. Recall that $\cM_{W^j}$ is the total occupation measure of $\widehat{W}^j$ and that $\cR_{W^j}$ is the range of $\widehat{W}^j$: 
\begin{equation}
\label{defavecj}
\cR_{W^j} = \left\{ \widehat{W}^j_t \; ;  \; t \in [0, \sigma_j ]  \right\} \quad {\rm and} \quad \langle \cM_{W^j} , f \rangle = \int_0^{\sigma_j } f ( \widehat{W}_t^j) \, dt \; , 
\end{equation}
where $\sigma_j$ stands for the total duration of $ \widehat{W}^j$.   
Recall the definition of $\cM^*_a$ from (\ref{Mstardef}). To simplify notation, we have assumed that these random variables are defined on an auxiliary probability space $(\Omega, \cG, \bP_0)$ and we also recall that under $\bP_0$, $\xi = (\xi_t, t \geq 0)$ is distributed as a $d$-dimensional Brownian motion starting at the origin. We first prove the following lemma. 

\begin{lemma}
\label{zeroone} There exists $\kappa_d  \in [0, \infty]$ that only relies on space dimension $d$ such that for any $a >0$, 
$$ \bP_0 {\rm -a.s.} \quad  \liminf_{r \rightarrow 0+} \; \frac{\cM^*_a (B(0, r))}{g(r)} = \kappa_d  \; .$$
\end{lemma}
\noi 
{\bf Proof:}  we first need to set some notation. For any $j \in \{1, \ldots, d \}$, $\xi^{(j)}$ stands for the $j$-th component process in the canonical basis of $\bR^d$; thus under $\bP_0$, the $\xi^{(j)}$'s,  are $d$ independent linear Brownian motions. For any $R \in [0, \infty)  $, we set 
\begin{equation}
\label{gammadef}
 \gamma (R)= \sup \left\{ t \geq 0 \; : \; \sqrt{ (\xi^{_{(1)}}_t)^2 + (\xi^{_{(2)}}_t)^2 + (\xi^{_{(3)}}_t)^2   } \, \leq \, R \;  \right\}  \; .
 \end{equation}
The process $\gamma=(\gamma (R), R \geq 0)$ is distributed as the three-dimensional Brownian escape process. Then, by a result due to Pitman \cite{Pit75}, $\gamma$ is a subordinator whose Laplace exponent is $\sqrt{2\lambda}$. Moreover, it enjoys  the following independence property: 
for any $R_1 < R_2 $, under the probability measure $\bP_0$, 
\begin{equation}
\label{escapeindep}
\gamma (R_1) \; , \; \gamma (R_2) -\gamma (R_1) \; {\rm and} \; \left( \xi_{ \gamma (R_2) +t } \, , \, t \geq 0 \right) \; {\rm are} \; {\rm independent}. 
\end{equation}
Let us fix $s >0$ and $r \in (0, 1)$. We define the following event 
$$ A (s,r)= \bigcap_{ \substack{ j\in \cJ^* \\ t_j > s}} \left\{  \cR_{W^j} \cap \overline{B}(0, r) = \emptyset   
\right\} \; .$$
We first claim that 
\begin{equation}
\label{localisationb}
 \forall s \in (0, \infty) \; , \quad \lim_{r \rightarrow 0} \bP_0 (\, A(s,r) \, )= 1 \; . 
\end{equation}
Indeed, recall notation $u_{0,r}$ from (\ref{uxrdef}); the definition of $\cN^*$ easily entails the following. 
\begin{eqnarray}
\label{separation} 
\bP_0 (\, A(s,r) \, ) & \geq & \bP_0 \left( \, A(s,r) \cap \{ s > \gamma (\sqrt{r})  \}  \right) \nonumber \\
& \geq & \bE_0 \left[ \un_{ \{ s > \gamma (\sqrt{r})  \}} e^{-4 \int_s^{\infty}  u_{0,r} (\xi_t) \, dt  } \right]  \nonumber  \\ 
& \geq & \bE_0 \left[ \un_{ \{ s >\gamma (\sqrt{r})  \}} e^{-4 \int_{\gamma (\sqrt{r})}^{\infty}  u_{0,r} (\xi_t) \, dt  } \right] . 
\end{eqnarray}
Next, we use (\ref{escapeindep}) to get 
$$ \bE_0 \left[ \un_{ \{ s > \gamma (\sqrt{r}) \}} e^{-4 \int_{\gamma (\sqrt{r})}^{\infty}  u_{0,r} (\xi_t) \, dt  } \right]  = \bP_0 ( s > \gamma (\sqrt{r}) \, )   \bE_0 \left[  e^{-4 \int_{\gamma (\sqrt{r})}^{\infty}  u_{0,r} (\xi_t) \, dt  } \right] . $$
Then, Jensen inequality, combined with Lemma \ref{expectgreen} with $b=r^{1/2}$, entails the following inequalities for any $r \in (0, 1/4)$. 
\begin{eqnarray*}
 \bE_0 \left[ e^{-4 \int_{\gamma (\sqrt{r})}^{\infty}  u_{0,r} (\xi_t) \, dt  } \right]  & \geq & \bE_0 \left[ e^{-4 
 \int_{0}^{\infty} \un_{ \{ \lVert \xi_t \lVert \geq \sqrt{r}  \}} u_{0,r} (\xi_t) \, dt  } \right] \\ 
& \geq & \exp \left( - 4  \bE_0 \left[
 \int_{0}^{\infty} \un_{ \{ \lVert \xi_t \lVert \geq \sqrt{r} \}} u_{0,r} (\xi_t) \, dt  \right] \right)  \\
 & \geq & \exp (- 4 C_7 r^{\frac{d-4}{2}} \, ) \; .  
\end{eqnarray*}
Therefore (\ref{separation}) implies 
$$ \bP_0( A(s,r)) \geq \bP_0 ( s > \gamma (\sqrt{r})  \, ) \cdot \exp (- 4 C_7 r^{\frac{d-4}{2}} \, ) \; , $$
which easily entails Claim (\ref{localisationb}).

Next, observe that $A(s, r) \subset A(s, r')$ if $r'<r$. Therefore, (\ref{localisationb}) and Borel-Cantelli lemma imply that for any fixed $s \in (0, \infty)$, $\bP_0 {\rm -a.s.}$ the event $A(s,r)$ is realised for all sufficiently small $r$. By definition of $A(s,r)$, it entails that for any $a>s>0$, 
\begin{equation}
\label{localisationc}
 \bP_0 {\rm -a.s.} \quad \liminf_{r \rightarrow 0}    \frac{\cM^*_a (B(0, r))}{g(r)} = \liminf_{r \rightarrow 0}    \frac{\cM^*_s(B(0, r))}{g(r)} \; .
 \end{equation}
Let us introduce the filtration $(\cG_s , s \geq 0)$ where $\cG_s$ is the sigma field generated by 
$\un_{[0,s]} (t)  \cN^* (dtdW) $ and completed with the $\bP_0$-negligible sets. Standard arguments on Poisson point measures combined with Blumenthal zero-one law for $\xi$ entail that the sigma field 
$\cG_{0+}= \bigcap_{s >0} \cG_s$ is trivial. This combined with (\ref{localisationc}) show there exists 
$\kappa_d \in [0, \infty]$ (that only relies on space dimension) such that  for any $a>0$, 
$\bP_0$-a.s. $ \kappa_p=  \liminf_{r \rightarrow 0}   g(r)^{-1} \cM^*_a (B(0, r))$, which is the desired result. \cqfd

\vspace{5mm}

 By (\ref{Palm}), the previous lemma entails that 
\begin{equation} 
\label{Partiel}
\bbN_0 {\rm -a.e.} \; {\rm for} \; \cM\!{\rm -almost} \; {\rm all} \; x,  \quad  \liminf_{r \rightarrow 0+} \frac{\cM 
(B(x, r))}{g(r)} = \kappa_d  \in [0, \infty] \; .
\end{equation}
Then to complete the proof of Theorem \ref{typilowdensM},  we only need to prove that $0 <\kappa_d  < \infty$, which is done in two steps.

\begin{lemma}
\label{upperbound} For any $d \geq 5$, $\kappa_d \leq 27/2$. 
\end{lemma}
\noi
{\bf Proof:} by Lemma \ref{zeroone}, we only need to prove that for any $a>0$, 
\begin{equation} 
\label{Palmifie}
\bP_0 \left( \liminf_{r \rightarrow 0+} g(r)^{-1} \cM^*_a (B(0, r)) \leq 27/2 \right) \; > \; 0 \; .
\end{equation}
To that end, we need to introduce the following notation: we recall the definition of $\gamma$ from (\ref{gammadef}) and we first set for any $r \in (0, \infty) $:
$$ S_r = \sum_{j \in \cJ^*} \un_{[ \, 0, \gamma (2r)\,  ]} (t_j) \,  \sigma_j  \; , $$
where $\sigma_j$ stands for the duration of $W^j$. Recall from (\ref{defavecj}) notation $\cM_{W^j}$ and $\cR_{W^j}$. We next define the following event: 
\begin{equation}
\label{erdefff}
 E_r =\bigcap_{ \substack{ j\in \cJ^* \\ t_j > \gamma (2r) }} \left\{  \cR_{W^j} \cap \overline{B}(0, r) = \emptyset   \right\} \; .
 \end{equation}
Recall that for any $j \in \cJ^*$, $\cM_{W^j} (B(0, r) ) \leq \sigma_j $. Consequently, 
\begin{equation}
\label{remarquefonda}
\bP_0{\rm -a.s.} \; {\rm on} \; E_r \; , \quad \cM^*_a (B(0,r)) \leq S_r \; .
\end{equation}
If we set for any 
any $n \geq 2$, 
$$r_n = \left( 1/\log n\right)^n  \quad {\rm and } \quad  V_n = \un_{ \{ \, S_{r_n} \leq  \frac{27}{2} \, g(r_n) \,    \} \cap E_{r_n}  }  \; , $$
then, (\ref{Palmifie}) is a consequence of the following: 
\begin{equation}
\label{discrete}
\bP_0 \left( \sum_{n \geq 2}  V_n = \infty \right) \; > \; 0 \; .
\end{equation}
Therefore, we only need to prove (\ref{discrete}). We proceed in three steps. 

\medskip

\noi
$\bullet$ {\it Step I:} we first claim there exists $C_{10} \in (0, \infty)$ that only relies on space dimension $d$ such that for any $r >0$: 
\begin{equation}
\label{localite}
\bP_0 (E_r) \geq  C_{10}\; .
\end{equation}
Indeed, from the definition of $\cN^*$, we get 
$$ \bP_0 (E_r)= \bE_0 \left[  e^{-4\int^\infty_{\gamma (2r) } u_{0,r} (\xi_t) \, dt  }\right] \; .$$
An easy argument combined with Jensen inequality entail the following. 
\begin{eqnarray*}
 \bP_0 (E_r) & \geq &  \bE_0 \left[  e^{-4\int_{0}^\infty \un_{ \{ \lVert \xi_t \rVert \geq 2r \}} u_{0,r} (\xi_t) \, dt  }\right] \\
& \leq & \exp \left( -4 \bE_0 \left[\int_{0}^\infty \un_{ \{ \lVert \xi_t \rVert \geq 2r \}} u_{0,r} (\xi_t) \, dt   \right] \right).
\end{eqnarray*}
Then, we apply Lemma \ref{expectgreen} with $b=2r$ to get $\bP_0 (E_r) \geq  \exp (-4C_7 2^{d-4}) := C_{10}$, which proves (\ref{localite}). 


\medskip

\noi
$\bullet$ {\it Step II:} we set $ L_n =  V_2 + \ldots + V_n$ and we claim that 
\begin{equation}
\label{suminfinite}
\lim_{n \rightarrow \infty}\bE_0 \left[  L_n \right] = \infty \; .
\end{equation}
To that end, we explicitely compute the distribution of the process $r \mapsto S_r$: since conditionally given $\xi$, $\cN^*$ is distributed as a Poisson point measure with intensity $4\, dt \, \bbN_{\xi(t)} (dW)$ and since $S_r$ only relies on $\xi$ via $\gamma (2r)$, (\ref{escapeindep}) easily implies that 
\begin{equation}
\label{indepexcape}
S_{r_1} \; , \; S_{r_2} -S_{r_1} \; {\rm and} \; \un_{E_{r_2}} {\rm are} \; {\rm independent}. 
\end{equation}
Next, recall that for any $x\in \bR^d$ and any $\lambda \geq0$, we have 
$$ \bN_x \left( 1-e^{-\lambda \sigma } \right) = N\left( 1-e^{-\lambda \sigma } \right)= \sqrt{\lambda/2} \; .$$
Thus, the exponential formula for Poisson point measures entails 
$$ \bE_0 \left[ e^{-\lambda (S_{r_2} -S_{r_1})}\right] = \bE_0 \left[ e^{-\sqrt{8\lambda}(\gamma (2r_2) -\gamma (2r_1) \,  )} \right] = e^{- (r_2-r_1) (128 \lambda)^{1/4} }\; .$$
This, combined with the independence property (\ref{indepexcape}), entails that $(S_r , r\geq0)$ is a stable subordinator with exponent $1/4$ and speed $128$. To prove (\ref{suminfinite}), we use the following estimate of the tail at $0+$ of $S_1$ that is due to Shorokhod \cite{Sko54} (see also Example 4.1 Jain and Pruitt \cite{JaPr87}).   
\begin{equation}
\label{tailestim}
\bP_0  \left( S_1 \leq x \right)  \; \sim_{x \rightarrow 0+} \; C_{11} x^{\frac{1}{6}} \exp \left( -(x/ C_{12})^{-\frac{1}{3}}\right) \; , 
\end{equation}
where, to simplify notation, we have set $C_{11}= (6\pi)^{-1/2}2^{7/6}$ and $C_{12}= 27/2$. 
Next, by (\ref{indepexcape}), the scaling property for $S$ and (\ref{localite}), we get 
\begin{eqnarray}
\label{ineqsplitage} 
\bE_0 \left[ V_n \right] & =& \bP_0( \, S_{r_n} \leq C_{12}\, g(r_n)\,  ) \cdot \bP_0 (E_{r_n}) \nonumber\\
&\geq &C_{10}  \bP_0 \left( \,  S_1 \leq C_{12} (\log \log 1/r_n)^{-3} \right)   \; .
\end{eqnarray}
Now by (\ref{tailestim}), we get 
\begin{equation}
\label{keyequival}
 \bP_0 \left(  S_1 \leq C_{12}(\log \log 1/r_n)^{-3} \right)  \, \sim_{n \rightarrow \infty} \,  \frac{ C_{11}C_{12}^{1/6} }{n\sqrt{\log n} \, \log \log n} \; , 
 \end{equation}
which easily implies  (\ref{suminfinite}).

\medskip

\noi
$\bullet$ {\it Step III:} we finally claim that there exists $C_{13} \in (0, \infty)$ that only depends on space dimension $d$, such that 
\begin{equation}
\label{claimdeux}
\forall \, 2 \leq k < \ell \; , \quad \bE_0 \left[ V_k V_\ell \right] \leq C_{13} \cdot  \bE_0 \left[ V_k \right] \bE_0 \left[  V_\ell \right] \; .
\end{equation}
Indeed, by (\ref{localite}), (\ref{indepexcape}), (\ref{ineqsplitage}) and the scaling property of $S$, the following inequalities hold true. 
\begin{eqnarray}
\label{ineqserie}
\bE_0 \left[ V_k V_\ell \right] & \leq & \bP_0 \left(  \,  S_{r_\ell} \leq C_{12}\, g(r_\ell)\; ; \; S_{r_k} -S_{r_\ell} \leq C_{12} \, g(r_k)\,  \right) \nonumber \\
& \leq &  \bP_0 \left(  \,  S_{r_\ell} \leq C_{12} \, g(r_\ell)\, \right) \cdot \bP_0 \left( S_{r_k} -S_{r_\ell} \leq C_{12} \, g(r_k) \right) \nonumber \\
& \leq & \frac{1}{C_{10}} \bE_0 \left[ V_\ell\right] \cdot  \bP_0 \left( S_1 \leq C_{12}  \, (1-r_{k}/r_{k+1})^{-4} (\log \log 1/r_k)^{-3} \,   \right) . 
\end{eqnarray}
An easy computation entails that 
$$   \bP_0  \left( S_1 \leq C_{12}\,  (1-r_{k}/r_{k+1})^{-4} (\log \log 1/r_k)^{-3} \,   \right) \sim_{k \rightarrow \infty} \,  \frac{ e^{4/3}C_{11}C_{12}^{1/6} }{k\sqrt{\log k} \, \log \log k} \; . $$
Then by (\ref{keyequival}) and (\ref{ineqsplitage}),  there exist $C_{14} \in (0, \infty) $ that only depends on space dimension $d$,  such that for any $k \geq 2$ 
$$  \bP_0 \left( S_1 \leq C_{12} \,  (1-r_{k}/r_{k+1})^{-4} (\log \log 1/r_k)^{-3} \,   \right) \leq C_{14} \,  \bE_0 
\left[ V_k \right] \; , $$
which entails (\ref{claimdeux}) by (\ref{ineqserie}) with $C_{13}= C_{14}/C_{10}$. 

\medskip

   Claim (\ref{suminfinite}) and Claim (\ref{claimdeux}) entail 
  $$ \limsup_{n \rightarrow \infty} \frac{ \bE_0 \left[ L_n^2 \right] }{ \bE_0 \left[ L_n \right]^2 } \leq C_{13}  
  \; , $$
and (\ref{discrete}) follows by Kochen-Stone's lemma, which completes the proof of the lemma. \cqfd

\vspace{5mm}

The following lemma completes the proof of Theorem \ref{typilowdensM}.  
\begin{lemma}
\label{lowerbound} For any $d \geq 5$, $\kappa_d \geq 2^{-10}$. 
\end{lemma}
\noi
{\bf Proof:} we directly work with $W$ under $\bbN_0$. Then, it is sufficient to prove that 
\begin{equation}
\label{directlift}
\bbN_0 {\rm -a.e.} \, {\rm for}  \; \cM\!{\rm -almost} \; {\rm  all} \;  x, \quad   \liminf_{r \rightarrow 0+} \; \, \frac{\cM (B(x, r))}{g(r)} \geq 2^{-10}\; .
\end{equation}
The proof of (\ref{directlift}) consists in lifting to $\widehat{W}$ estimates from Lemma \ref{typidensBtree} by using the fact that conditionally given $H$, $\widehat{W}$ is a Gaussian process. Recall notation $d_H$. For any $r, R >0$   and for any $t \in [0, \sigma]$, we set  
$$ {\bf a}(t,r)= \int_0^\sigma \un_{ \{ d_H (s,t) \leq r \}} \, ds \quad {\rm and} \quad  {\bf b}(t,r, R)= \int_0^\sigma 
\un_{ \{ d_H (s,t) \leq r \} \cap \{  \lVert \widehat{W}_s -\widehat{W}_t \rVert  \geq R \}} \, ds \; .$$
Then for any $t \in [0, \sigma]$, we first notice the following. 
\begin{equation}
\label{ineqdens}
 {\bf a}(t,r) \leq {\bf b} (t,r, R) + \cM (B(\widehat{W}_t , R)) \; .
\end{equation}
Recall that $\widehat{W}$, conditionnaly given $H$, is distributed as a centered  Gaussian process whose covariance is specified by (\ref{covar}). Consequently, 
\begin{equation}
\label{condexpec}
 N(dH) {\rm -a.e.} \; \forall t \in [0, \sigma ] \; , \; \,  Q^H_0 \left[ {\bf  b} (t,r, R) \right] \leq {\bf a}(t,r)  \int_{ \bR^d \backslash B(0, R/\sqrt{r})} \! \! \! \! \! \! \! \! \!\! \! \! \! \! \! \! \! \! (2\pi)^{-d/2} e^{-\lVert x \rVert^2/2 } dx \; .
\end{equation}
Next, for any integer $n \geq 2$, we set $R_n= 2^{-n}$ and $r_n= \frac{1}{4} R_n^2 (\log \log 1/R_n)^{-1}$. An elementary argument implies that for any $n \geq 2$, 
$$ \int_{ \bR^d \backslash B(0, R_n/ \sqrt{r_n})} \! \! \! \! \! \! \! \! \!\! \! \! \! \! \! \! \! \! (2\pi)^{-d/2} e^{-\lVert x \rVert^2/2 } dx \; \,  \leq \; \, C_{16} n^{-3/2} \; , $$ 
where $C_{16}$ is a positive number that only depends on space dimension $d$ (note that the power $3/2$ in the previous inequality is not optimal). Therefore, we get 
$$  N(dH) {\rm -a.e.} \; \forall t \in [0, \sigma ] \; , \quad  Q^H_0 \left[ \sum_{n\geq 2} \frac{{\bf b}(t,r_n, R_n)}{{\bf a}(t,r_n)} \right] 
 < \infty \; .$$
 Then, by Fubini, 
 $$ N(dH) {\rm -a.e.}  \; \quad Q^H_0 \left[ \int_0^{\sigma} \un_{ \left\{ \limsup_{n \rightarrow \infty}  \frac{{\bf b}(t,r_n, R_n)}{{\bf a} (t,r_n)} \; \, > \; \, 0  \right\}} \, dt  \right]= 0 \; , $$
which implies that 
\begin{equation}
\label{almostneglig}
 \bbN_0-{\rm a.e.} \; {\rm for} \; \ell \! {\rm -almost} \; {\rm all}\; t \in [0, \sigma ]\, , \;   \quad 
\lim_{n \rightarrow \infty}  \frac{{\bf b}(t,r_n, R_n)}{{\bf a} (t,r_n)}  = 0  
\end{equation}
($\ell$ stands here for the Lebesgue measure on the real line). Recall notation $k(r)$ from (\ref{defkk}); 
(\ref{almostneglig}) combined with
(\ref{ineqdens}) entails 
$$  \bbN_0-{\rm a.e.} \; {\rm for} \; \ell \! {\rm -almost} \; {\rm all}\; t \in [0, \sigma ],  \; \, 
\liminf_{n \rightarrow \infty}  \frac{{\bf a}(t,r_n)}{k(r_n)}  \leq \liminf_{n \rightarrow \infty}  \frac{\cM (B(\widehat{W}_t , R_n))}{k(r_n)}  \; .$$
Since $k(r_n) \sim_{n \rightarrow \infty}2^{-4} g(R_n) $, Lemma \ref{typidensBtree} entails that 
$$  \bbN_0-{\rm a.e.} \; {\rm for} \; \ell \! {\rm -almost} \; {\rm all}\; t \in [0, \sigma ], \; \,  
\liminf_{n \rightarrow \infty}  \frac{\cM (B(\widehat{W}_t ,  2^{-n}))}{g(2^{-n})} \geq 2^{-6} , $$
and an easy argument completes the proof of the lemma. \cqfd

\subsection{Proof of Theorem \ref{packingsnake}. }

We first introduce a specific decomposition of $\bR^d$ into dyadic cubes. We adopt the following notation: we denote by $\lfloor \, \cdot \rfloor$ the integer part application and we write $\log_2$ for the logarithm in base $2$; we fix $d \geq 5$ and we set 
$$ p := \lfloor \log_2 (4\sqrt{d} ) \rfloor \; , $$ 
so that $2^p > 2\sqrt{d}$. To simplify notation, we set $\cD_n = 2^{-n-p}\bZ^d $, for any $n \geq 0$. 
For any $y = (y_1, \ldots , y_d) $ in $\cD_n$, we also set 
$$ D_n (y)= \prod_{j= 1}^d \, [ \, y_j -\frac{_1}{^2} 2^{-n} \,  ;  \, y_j+\frac{_1}{^2} 2^{-n} \, ) \;  \quad {\rm and} \quad 
D^\bullet_n (y)= \prod_{j= 1}^d\,  [ \, y_j - \frac{_1}{^2} 2^{-n-p} \, ; \, y_j +\frac{_1}{^2} 2^{-n-p} \, )  .$$
It is easy to check the following properties. 
\begin{itemize}
\item{Prop(1).} If $y,y'$ are distinct points in $\cD_n$, then $D^\bullet_n (y) \cap D^\bullet_n (y') = \emptyset$. 

\medskip

\item{Prop(2).} Let $y \in \cD_n$. Then, we have 
$$ D^\bullet_n (y) \,  \subset \, \overline{B} ( y \, , \, \frac{_1}{^2} 2^{-n-p}\sqrt{d} \, ) \, \subset  \, 
\overline{B} ( y \, , \, 2^{-n-p}\sqrt{d} \, ) \, \subset \, D_n (y) \; .$$
\end{itemize}
 For any $ r < (2d)^{-1}$, we set $ n(r)= \lfloor \log_2 ( r^{-1} (1+2^{-p})\sqrt{d}  )\rfloor $, so that the following inequalities hold: 
\begin{equation}
\label{coincr}
  \frac{_1}{^2} (1+2^{-p})\sqrt{d} \cdot 2^{-n(r)} < r \leq (1+2^{-p})\sqrt{d} \cdot 2^{-n(r)} \; .
\end{equation}
Next, for any $x = (x_1, \ldots , x_d)\in \bR^d$ and for $j \in \{ 1, \ldots , d \}$, we set 
$$ y_j = 2^{-n(r)-p} \lfloor  x_j 2^{n(r)+p} +  \frac{_1}{^2} \rfloor  \; .$$
Therefore, $y= (y_1, \ldots , y_d) \in \cD_{n(r)}$ and we easily check the following: 
\begin{itemize}
\item{Prop(3).} The point $x$ belongs to  $D_{n(r)}^\bullet (y) $ and $D_{n(r)} (y) \,  \subset\,  B(x,r) $. 
\end{itemize}

\medskip

Recall that we work under $\bbN_0$ and recall $C_3$ and $\kappa_0$ from Lemma \ref{endestimate}. 
We set $\kappa_1 = \min ( \kappa_0 /2 \, , \, C_3/8)$ and we choose $\kappa_2 >0$ such that for 
\begin{equation}
\label{tutune}
\forall n \geq 7 \; : \quad  \kappa_2 g(2^{-n}) \leq \kappa_1 g( \frac{_1}{^2} 2^{-n-p} \sqrt{d} ) \; . 
\end{equation}
We then fix $A>100$ and  for any $n $ such that $2^{-n} \leq 1/(2A)$, we set 
$$ U_n (A) = \sum_{\substack{ y \in \cD_n \\ 1/A \leq \lVert y \rVert \leq A}} g( \sqrt{d} (1+2^{-p}) 2^{-n} ) \un_{ \{ \,  \cM (D_n(y) ) \leq \kappa_2  g(2^{-n}) \,  \} \cap \{  \, \cR \cap D_n^\bullet (y) \neq \emptyset \, \} } \;. $$
We first prove the following lemma. 
\begin{lemma}
\label{controlvariable} For any $A>100$, $\bbN_0$-a.e. 
\begin{equation}
\label{Ulimit}
 \lim_{N \rightarrow \infty} \sum_{n \geq N} U_n (A) = 0 \; .
 \end{equation}
\end{lemma}
\noi
{\bf Proof:} we fix $n$ such that $2^{-n} \leq 1/(2A)$ and we fix $y \in \cD_n $ such that  $1/A\leq  \lVert y \rVert \leq A$. By Prop(2), we get 
\begin{eqnarray*}
 \bbN_0 \left(   \cM (D_n(y) ) \leq \kappa_2 g(2^{-n}) \; ; \;  \cR \cap D_n^\bullet (y) \neq \emptyset  \right)  & &
  \\
 & &   \hspace{-60mm}  \leq \bbN_0 \left(   \cM (  \overline{B} ( y \, , \, 2^{-n-p}\sqrt{d} \, ) ) \leq \kappa_1 
 g( \frac{_1}{^2} 2^{-n-p} \sqrt{d}) \; ; \;  \cR \cap \overline{B} ( y \, , \, \frac{_1}{^2} 2^{-n-p}\sqrt{d} \, )   \neq \emptyset  \right) . 
\end{eqnarray*}
We next apply Lemma \ref{endestimate} with $x= y$ and $r= \frac{_1}{^2} 2^{-n-p}\sqrt{d}$ to prove there exists  $C_{17} \in (0, \infty)$, that only depends on space dimension $d$, such that: 
$$   \bbN_0 \left( \,   \cM (D_n(y) ) \leq \kappa_1 g(2^{-n}) \; ; \;  \cR \cap D_n^\bullet (y) \neq \emptyset \,   \right)  
\leq C_{17} (2^{-n-p})^{d-4} \lVert y \rVert^{2-d} n^{-2} \; .$$
Then, note there exists $C_{18} \in (0, \infty) $, that only depends on space dimension $d$, such that for all sufficiently large $n$, 
$$ g( \sqrt{d} (1+2^{-p})  2^{-n} ) \leq C_{18} (2^{-n-p})^{4}  \; , $$
which entails the following. 
\begin{eqnarray*}
g( \sqrt{d} (1+2^{-p})2^{-n} ) \cdot  \bbN_0 \left(  \,  \cM (D_n(y) ) \leq \kappa_1 g(2^{-n}) \, ; \,  \cR \cap D_n^\bullet (y) \neq \emptyset \,  \right)  & &\\
& & \hspace{-40mm}  \leq C_{19} (2^{-n-p})^{d} \lVert y \rVert^{2-d} n^{-2} \; , 
\end{eqnarray*} 
where $C_{19}= C_{17} C_{18}$. Elementary arguments entail the following inequalities.  
\begin{eqnarray*}
 \bbN_0 \left(  U_n (A) \right)& \leq & C_{19}\,  n^{-2} \sum_{\substack{ y \in \cD_n \\ 1/A \leq \lVert y \rVert \leq A} } (2^{-n-p})^{d} \lVert y \rVert^{2-d} \\
 & \leq & C_{20}\,  n^{-2} \int \un_{ \{1/A\leq  \lVert x \rVert \leq A \}}  \lVert x \rVert^{2-d} dx \\
 & \leq & C_{21} \, n^{-2} \int_{1/A}^A \rho \, d\rho \\
& \leq &  C_{21} \, A^2 n^{-2} , 
\end{eqnarray*}
where $C_{20},  C_{21} \in (0, \infty) $ only depends on $d$. Therefore,  
$ \bbN_0 \left( \sum_{n \geq N} U_n (A) \right) < \infty  $, which easily completes the proof of the lemma.  \cqfd

\medskip

  We next prove the following lemma. 
\begin{lemma}
\label{badpoint} Assume that $d \geq 5$. Then, $\bbN_0$-a.e. we have 
\begin{equation}
\label{badzero}
 \cP_g \left(  \left\{ x \in \cR \; :\; \liminf_{r \rightarrow \infty} g(r)^{-1} \cM (B(x,r)) \, \neq \kappa_d  \, \right\} \right) = 0 \; .
 \end{equation}
\end{lemma}
\noi
{\bf Proof:} we fix $A >100$. By Theorem \ref{typilowdensM} and Lemma \ref{controlvariable} there exists a Borel subset $\cW_A$ of $\cW$ such that $\bbN_0 (\cW \backslash \cW_A)= 0$ and such that on 
$\cW_A$,  (\ref{typidensMM}) and (\ref{Ulimit}) hold true. We argue for a fixed $W \in \cW_A$.

   Let $B $ be any  Borel subset of  $ \cR \cap  \{ x \in \bR^d \, : \, 1/A \leq \lVert x \rVert   \leq A  \} $. 
Let $\varepsilon >0$ and let $\overline{B} (x_1, r_1), \ldots \overline{B} (x_k, r_k)$ be any closed $\varepsilon$-packing of $B \cap \cR$.  Let $C_{22} $ be a positive real number to be specified later. 
First observe that 
\begin{eqnarray}
\label{splitineq}
\sum_{i=1}^k  g(r_i)  &= & \sum_{i=1}^k  g(r_i)  \un_{ \{  \cM ( B( x_i, r_i) >C_{22} \, g(r_i)  \}} +  \sum_{i=1}^k  g(r_i)  \un_{ \{  \cM ( B( x_i, r_i) \leq C_{22} \, g(r_i)  \}} \nonumber \\
& \leq & C_{22}^{-1} \cM \left( \, B^{(\varepsilon)} \, \right) +  \sum_{i=1}^k  g(r_i)  \un_{ \{  \cM ( B( x_i, r_i) \leq C_{22} \,  g(r_i)  \}} \; , 
\end{eqnarray}
where for any bounded subset $B$ of $\bR^d $ we have set $B^{(\varepsilon)}= \{ x \in \bR^d \; : \; {\rm dist} (x, B) \leq \varepsilon \} $. 

\medskip

  Next, fix $1 \leq j \leq k$; recall notation $n(r_i)$ from (\ref{coincr}) and denote by $y_i$ the point of $\cD_{n(r_i)}$ corresponding to $x_i$ such that Prop(3) holds true. Therefore, by (\ref{coincr}), we have 
\begin{eqnarray*}
 \cM ( B( x_i, r_i)\, )  \leq C_{22} \, g(r_i) \quad {\rm and} \quad  \;   x_i \in B \cap \cR \quad  \Longrightarrow  & & \\
& & \hspace{-70mm}  \cM (D_{n(r_i)} (y_i) )\leq C_{22} \, g( (1+2^{-p})\sqrt{d} 2^{-n(r_i) }) \quad {\rm and}  \quad \cR \cap D^\bullet_{n(r_i)} (y_i) \neq \emptyset .
\end{eqnarray*}
We now choose $C_{22}$ in order to have: $ C_{22} g( (1+2^{-p}) \sqrt{d} 2^{-n(r) }) \leq \kappa_2 g( 2^{-n(r) }) $ for all sufficiently small $r\in (0, 1)$. Thus, we get 
$$ \sum_{i=1}^k  g(r_i)  \un_{ \{  \cM ( B( x_i, r_i) \leq C_{22} g(r_i)  \}} \;  \leq \sum_{n\, : \,  2^{-n}  \leq C_{23} \, \varepsilon } U_n (A) \; , $$
where $C_{23}= 2((1+2^{-p}) \sqrt{d})^{-1}$. Since $W $ belongs to $\cW_A$ where (\ref{Ulimit}) holds, this inequality combined (\ref{splitineq}) with implies the following. 
\begin{equation}
\label{abscont}
\cP_g \left( B \cap \cR\right) \leq \cP^*_g \left( B \cap \cR\right) \leq C_{22}^{-1}\;  \cM \left(\,  \bigcap_{\varepsilon >0}B^{(\varepsilon)} \right) \; .
\end{equation}
We next applies (\ref{abscont}) with  $B=B_A$ given by 
$$B_A = \left\{ x \in \cR \; : \;  1/A \leq \lVert x \rVert \leq A \; {\rm and} \; \liminf_{r \rightarrow 0} g(r)^{-1} \cM (B(x,r)) \neq \kappa_d    \right\} \; .$$
This shows that  $\cP_g (B_A)< \infty$. Suppose now that $\cP(B_A) >0$, then by (\ref{innerreg}), there exists a compact subset 
$K $ of $B_A$ such that 
$\cP_g (K) > 0$. Since $K$ is compact then $K=\bigcap_{\varepsilon >0} K^{(\varepsilon)}$; now, since 
$K$ is a subset of $B_A$ and since $W\in \cW_A$ where (\ref{typidensMM}) holds true, we then get $\cM (K) =0$; by applying (\ref{abscont}) with $B=K$, we obtain $\cP_g (K)=0$, which rises a contradiction. Thus, we have proved that $\bbN_0$-a.e. $ \cP_g \left( B_A \right) \, = \, 0 $, 
which easily entails the lemma by letting $A$ go to $\infty$, since $\cP_g (\{ 0 \} ) = 0$. \cqfd 

\medskip

We now complete the proof of Theorem \ref{packingsnake}: by Theorem \ref{typilowdensM} and Lemma \ref{badpoint} there exists a Borel subset $\cW^*$ of $\cW$ such that $\bbN_0 (\cW \backslash \cW^*)= 0$ and such that (\ref{typidensMM}) and (\ref{badzero}) hold true on $\cW^*$. We fix $W \in \cW^*$ and we set 
$$ {\rm Good} = \left\{  x \in \cR \; : \;   \liminf_{r \rightarrow 0} g(r)^{-1} \cM (B(x,r)) = \kappa_d   \right\} \quad {\rm and } \quad {\rm Bad}= \cR \backslash {\rm Good}. $$
Let $B$ be any Borel subset of $\bR^d$. By (\ref{typidensMM}) and (\ref{badzero}), we have 
$$ \cM ( B\cap {\rm Bad}) =\cP_g (B\cap \cR \cap {\rm Bad}) = 0 \; .$$
Then, we apply Lemma \ref{equalitycomp} to ${\rm Good}\cap B$ and we get 
$$ \cM ( B \cap {\rm Good})= \kappa_d  \cdot \cP_g ( B \cap \cR \cap {\rm Good} ) \; .  $$
Therefore $\cM (B)= \kappa_d  \cdot \cP_g (B \cap \cR)$, which completes the proof of Theorem \ref{packingsnake}. \cqfd


\vspace{5mm}

\section{Proof of Theorem \ref{mainth}.}

 We derive Theorem \ref{mainth} from Theorem \ref{packingsnake}.  To that end, we first need an upper bound of the upper box-counting dimension of $\cR$ under $\bN_x$ . Let us briefly recall the definition of  the upper box-counting dimension: let $K$ be a compact subset of $\bR^d$; for any $\varepsilon >0$, we denote by $\cN (K, \varepsilon) $ the minimal number of balls with radius less than $\varepsilon$ that are necessary to cover $K$. The upper-box counting dimension of $K$ is then given by 
$$ \overline{{\rm dim}}_{{\rm Box}} (K)= \limsup_{ \varepsilon \rightarrow 0} \frac{\log \cN (K, \varepsilon) }{ \log (1/\varepsilon )} \; .$$
Let us fix $x \in \bR^d$. It is easy to prove any for any $ q \in (0, 1/4)$, $ \bN_x$-a.e. the endpoint process $(\widehat{W}_s, s\in [0, \zeta])$ is $q$-H\"older continuous (see Le Gall \cite{LG99} for a simple proof). This  implies that 
\begin{equation} 
\label{upperboxrange}
 \bN_x\! -{\rm a.e.} \quad \overline{{\rm dim}}_{{\rm Box}} (\cR) \leq 4 \; .
\end{equation}
(Actually, the Hausdorff, packing, upper and lower box-counting dimensions of $\cR$ are equal to $4$.) 
We prove the following lemma. 
\begin{lemma}
\label{capacity} Let $d\geq 5$ and let $x \in \bR^d$. For any compact set $K$ such that $\overline{{\rm dim}}_{{\rm Box}} (K) \leq 4$, we $\bN_x $-a.e. have $\cM (K)= 0$. 
\end{lemma}
\noi
{\bf Proof:} let us first assume that $x \notin K$ and set $k= \inf_{y \in K} \lVert x-y \rVert > 0$. For any 
$\varepsilon  \in (0, k/2)$, there exists $n_\varepsilon :=\cN (K, \varepsilon)$ balls denoted by 
 $B(x^\varepsilon_1 , \varepsilon )$, ..., $B( x^\varepsilon_{n_\varepsilon} , \varepsilon)$ that cover $K$. 
Then, (\ref{occupequa}) combined with standard estimates of $d$-dimensional Green function entail the following inequalities.
 \begin{eqnarray*}
 \bN_x ( \cM (K)) & \leq &  \sum_{i=1}^{n_\varepsilon} \bN_x \left( \cM(    B (x^\varepsilon_i , \varepsilon ) \right) \\
 & \leq &  \sum_{i=1}^{n_\varepsilon}  \int_0^\infty  \bP_x  \left[  \xi_t \in B (x^\varepsilon_i , \varepsilon ) 
    \right]  \, dt \\
 & \leq & C_{24} \,  \sum_{i=1}^{n_\varepsilon} 
 \int_{B (x^\varepsilon_i , \varepsilon )} \lVert x- y \rVert^{2-d} dy \\
 & \leq &  C_{25}  \, k^{2-d} \cdot  \varepsilon^{d}   n_\varepsilon ,  
\end{eqnarray*}
 where $C_{24}, C_{25} \in (0, \infty) $ only depend on $d$. Since we assume that 
 $d  > 4 \geq  \overline{{\rm dim}}_{{\rm Box}} (K)$, the previous inequality implies that $\bN_x (\cM (K))= 0$.

 Let us now consider the general  case: for any $r >0$, the previous case applies to the compact set $K'= K \backslash B(x,r)$ and we get:  
 $\bbN_x$-a.e. 
 $$ \cM (K) = \cM ( K \cap B(x, r)) + \cM (K \backslash B(x,r)) \leq \cM (B(x, r)) \; , $$ 
 and the proof of the claim is completed by letting $r$ go to $0$, since $\cM$ is diffuse. \cqfd
 
 \vspace{5mm}

The end of the proof of Theorem \ref{mainth} is a simple adaptation of Le Gall's argument in \cite{LGRange} (see p. 312-313) to the packing measure: first observe that we only need to consider the $\beta=1$ case for if we replace $Z$ by $c\cdot Z$, then ${\bf M}$ is replaced by $c\cdot {\bf M}$ but 
${\bf R}$ is unchanged.

  Let us consider the $\beta = 2$ case. Theorem \ref{packingsnake}  and Lemma \ref{capacity} imply that for any compact set $K$ such that $\overline{{\rm dim}}_{{\rm Box}} (K) \leq 4 $, and for any $x \in \bR^d$, we $\bN_x$-a.e. have: 
 \begin{equation}
 \label{packingpolar}
 \cP_{g} (K \cap \cR ) = 0 \; .
 \end{equation}
 Recall the connection (\ref{connection}) between ${\bf R}$, ${\bf M}$ and the excursions $W^j$, $j \in \cJ$, of the Brownian snake. An easy argument on Poisson point processes combined with (\ref{upperboxrange}) and  (\ref{packingpolar}) implies that almost surely $\cP_{g} \left( \cR_{W^{j}} \cap  \cR_{W^{i}}\right) =0$ for any $i\neq j $ in $\cJ$. Then, 
(\ref{connection}) entails  
$$ \cP_{g} \left( \, \cdot \cap {\bf R} \, \right) = \sum_{j \in \cJ } \cP_{g} \left( \, \cdot \cap \cR_{W^j}   \right) \; .$$
Theorem \ref{packingsnake} and (\ref{connection}) thus imply
$$ \kappa_d \cdot \cP_{g} \left( \, \cdot \cap {\bf R}  \right)  =  \sum_{j \in \cJ } \kappa_d  \cdot  \cP_{g} \left(\,  \cdot \cap \cR_{ W^{j} } \right)  =\sum_{j \in \cJ }  \cM_{W^{j} } = {\bf M}\; , $$
which is the desired result.  \cqfd 

%


%

%


%



%
\end{document}